\documentclass[twoside,a4paper,reqno,intlimits]{amsart}%
\usepackage[active]{srcltx}
\usepackage{color,esint}
\usepackage{amsthm}
\usepackage{amssymb}
\usepackage{amsfonts}

\usepackage{enumitem}
\setenumerate{label={\rm (\alph{*})}}

\usepackage[rose,frak,operator]{paper_diening}

\usepackage{graphicx}
\usepackage[%
   colorlinks=true,
   urlcolor=rltblue,       
   filecolor=rltgreen,     
   citecolor=rltgreen,     
   linkcolor=rltred,       
   bookmarks=true,
   unicode,
   plainpages=false,
   ]{hyperref}
\usepackage{xcolor}
\definecolor{rltred}{rgb}{0.75,0,0}
\definecolor{rltgreen}{rgb}{0,0.5,0}
\definecolor{rltblue}{rgb}{0,0,0.75}

\definecolor{egreen}{rgb}{0,0.6,0}
\newcommand\farb[1]{\color{black}#1 \color{black}}

\newtheorem{theorem}[equation]{Theorem}
\newtheorem{lemma}[equation]{Lemma}
\newtheorem{proposition}[equation]{Proposition}
\newtheorem{cor}[equation]{Corollary}
\newtheorem{assumption}[equation]{Assumption}
\newtheorem{remark}[equation]{Remark}

\numberwithin{equation}{section}

\newcommand{\meantmp}[2]{#1\langle{#2}#1\rangle}
\newcommand{\mean}[1]{\meantmp{}{#1}}

\newcommand{\Pidiv}{\ensuremath{\Pi^{\divergence}_h}}

\newcommand{\para}{{\delta}}
\newcommand{\PiY}{\Pi^Y_h}

\newcommand{\bfumh}{\mbox{\bf u}^m_h}

\newcommand{\dx}{{\rm d\bx}}

\newcommand{\bfuMh}{\mbox{\bf u}_h^m}
\newcommand{\bfeMh}{{\bf e}_h^{m}}
\newcommand{\bfeMsh}{{\bf e}_h^{m-1}}
\newcommand{\bfuMsh}{{\bf u}_h^{m-1}}
\newcommand{\piMh}{q _h^m}
\newcommand{\bfvM}{\mbox{\bf v}_h}

\newcommand{\Int}{\dashint\limits_{\!\!I_m}}
\newcommand{\ksum}{{\kappa\sum_{m=0}^M}}
\newcommand{\ksumo}{{\kappa\sum_{m=1}^M}}

\begin{document}
\title[Optimal error estimate for generalized Newtonian fluids] {Optimal error estimate
  for a space-time discretization for incompressible generalized Newtonian fluids: The
  Dirichlet problem}
\author{Luigi C.\ Berselli}
\address{Dipartimento di Matematica, Universit{\`a} di Pisa, Via F.~Buonarroti 1/c,
  I-56127 Pisa, ITALY.}  \email{luigi.carlo.berselli@unipi.it}

\author{Michael R\r u\v zi\v cka{}} 
\address{Institute of Applied Mathematics, Albert-Ludwigs-University Freiburg,
  Ernst-Zermelo-Str.~1, D-79104 Freiburg, GERMANY.}  \email{
  rose@mathematik.uni-freiburg.de}
\dedicatory{Dedicated to Hideo Kozono on the occasion of his 60th birthday}

\begin{abstract}
  In this paper we prove optimal error estimates for {solutions with natural regularity}
  of the equations describing the unsteady motion of incompressible shear-thinning fluids. We consider a full
  space-time semi-implicit scheme for the discretization.  The main novelty, with respect
  to previous results, is that we obtain the estimates directly without introducing
  intermediate semi-discrete problems\farb{, which enables the treatment of homogeneous
    Dirichlet boundary conditions.}
\end{abstract}
\date{\small \today}
\maketitle
\section{Introduction}
In this paper we study a space-time discretization of the unsteady
system describing the motion of homogeneous (for simplicity the
density $\rho$ is set equal to $1$), incompressible shear-thinning
fluids under homogeneous Dirichlet boundary conditions. We prove
optimal error estimates (cf.~Section~\ref{sec:main}) for solutions
possessing a \textit{natural} regularity, extending the results in~\cite{bdr-3-2} to the
case of homogeneous Dirichlet boundary 
conditions. Our method differs from most previous investigations in as
much as we use no intermediate semi-discrete problems to prove our
result. We restrict ourselves to the three-dimensional setting,
however, all results can be easily adapted to the general setting in
$d$-dimensions.

More precisely, we consider for a bounded polyhedral domain
$\Omega \subset \setR^3$ and a finite time interval $I:=(0,T)$, for
some given $T>0$, the system
\begin{equation}
  \label{eq:pfluid}
  \begin{aligned}
    \partial_t\bfu-\divergence \bfS(\bfD\bfu)+ [\nabla \bfu] 
    \bfu+\nabla q&= \bff\qquad&&\text{in } I\times \Omega,
    \\
    \divergence\bfu&=0\qquad&&\text{in } I\times \Omega, 
    \\
    \bfu(0)&=\bfu_0\qquad&&\text{in }  \Omega,
  \end{aligned}
\end{equation}
where the vector field $\bfu=(u_1,u_2,u_3)^\top$ is the velocity, the
scalar $q$ is the kinematic pressure, the vector
$\bff=(f_1,f_2,f_3)^\top$ is the external body force and  $\bfu_0$ is the
initial velocity. We assume that the extra stress tensor
$\bS$ has $(p,\delta)$-structure \farb{for some $p\in(1,2]$, and
  ${\delta\in[0,\infty)}$} (cf.~Section~\ref{sec:stress_tensor}). For
the convective term we use the notation
$([\nabla \bfu]\bfu)_i = \sum_{j=1}^3 u_j \partial_j u_i$,
${i=1,2,3}$, while $\bfD \bfu := \tfrac 12 (\nabla\bfu+\nabla
\bfu^\top )$ denotes the symmetric part of the gradient $\nabla \bfu$.
For smooth enough solutions $(\bu,q)$ the variational formulation
of~\eqref{eq:pfluid} reads
\begin{align}    
  \label{eq:cont-var}
  \begin{aligned}
    \hskp{ {\partial _t\bfu}(t)
    }{\bfv}+\hskp{\bS(\bD\bfu(t))}{\bD\bfv}+b(\bu(t),\bu(t),\bv)-\hskp{q(t)}{\divo
      \bv}&=\hskp{\bff(t)}{\bfv}\,,
    \\
    \hskp{\divo \bu(t)}{\eta}&=0\,,
    \\
    \hskp{\bu(0)}{\bv}&=\hskp{\bu_0}{\bv}\,,
  \end{aligned}
\end{align}
for all $\bv \in
V:=W^{1,p}_{0}(\Omega)^3$, $\eta \in Q:=L_0^{p'}(\Omega)$
and almost every $t \in I$. We used the notation 
\begin{equation}\label{def:b}
 b(\bfu,\bfv,\bfw):=
  \frac{1}{2}\big[\hskp{[\nabla\bfv]\bfu}{\bfw}-\hskp{[\nabla\bfw]\bfu}{\bfv}\big],
\end{equation}
for the convective term to have a stable space-discretization.  Note
that $b(\cdot, \cdot,\cdot)$ is skew-symmetric with respect to the
last two arguments, i.e., $b(\bfu,\bfv,\bfw)\!=\!-b(\bfu,\bfw,\bfv)$
and that for solenoidal 
functions it holds $b(\bfu,\bfv,\bfw)= \hskp{[\nabla\bfv]\bfu}{\bfw}$.
We perform an error analysis for the semi-implicit space-time
discretization, which for given $h>0$ and $M \in \setN$ reads: for
$\bu_h^0 :=\Pidiv \bu_0$ find $(\bfuMh,\piMh) \in V_h\times Q_h$, $ m=1,\ldots,M$, such
that for all $\bfvM \in V_h$, $\eta _h \in Q_h$ holds 
\begin{align}  
  \label{eq:fdiscr-var}
  \begin{aligned}
    (d_t \bfuMh,\bfvM)+\hskp{\bS(\bD\bfuMh)}{\bfv_h}
    +b(\bfuMsh,\bfuMh,\bfvM)-\hskp{\piMh }{\divo
      \bfvM}&=(\bff(t_m),\bfvM),
    \\
    \hskp{\divo \bfuMh}{\eta_h}&=0\,,
  \end{aligned}
\end{align}
where $d_t \bu_h^m:=\kappa^{-1} ({\bu_h^m-\bu_h^{m-1}})$ is the backward
difference quotient with $\kappa :=\frac TM$. Here  $V_h\subset V$, $Q_h\subset Q$ are
appropriate stable finite element spaces with mesh size
$h>0$. The precise setup can be found in Section~\ref{sec:disc}.
\section{Preliminaries and main results}
\label{sec:preliminaries}
In this section we introduce the notation, the setup, and recall some
technical results which will be needed in the proof of the main result.
\subsection{Function spaces}
We use $c, C$ to denote generic constants, which may change from line
to line, but are not depending on the crucial quantities. 
We write $f\sim g$ if and
only if there exist constants $c,C>0$ such that
$c\, f \le g\le C\, f$.

We will use the customary Lebesgue spaces
$(L^p(\Omega), \norm{\,.\,}_p)$ and Sobolev spaces
$(W^{k,p}(\Omega), \norm{\,.\,}_{k,p})$, $k \in \setN$. We do not
distinguish between scalar, vector-valued or tensor-valued function
spaces in the notation if there is no danger of confusion. However, we
denote vector-valued functions by small boldfaced letters and
tensor-valued functions by capital boldfaced letters.  If a norm
is considered on a set $M$ different than $\Omega$ we indicate this
in the respective norms as $\norm{\,.\,}_{p,M},
\norm{\,.\,}_{k,p,M}$. We equip $W^{1,p}_0(\Omega)$ with the gradient
norm $\norm{\nabla \,.\,}_p$.  We denote by $\abs{M}$ the
$3$-dimensional Lebesgue measure of a measurable set $M$. The mean
value of a locally integrable function $f$ over a measurable set
$M \subset \Omega$ is denoted by
$\mean{f}_M:= \dashint\limits_M f \, dx =\frac 1 {|M|}\int\limits_M f \,
dx$. By $L^p_0(\Omega)$ we denote the space of functions
$f \in L^p(\Omega)$ with $\mean{f}_\Omega=0$.  Moreover, we use the
notation $({f},{g}):=\int\limits_\Omega f g\, dx$, whenever the
right-hand side is well-defined.

We use the following notation 
\begin{gather*}
  X :=\big(W^{1,p}(\Omega)\big)^3\,, \qquad  V
  :=\big(W^{1,p}_0(\Omega)\big)^3\,, \qquad V_{\divo}:=\bigset{\bv \in
  V \fdg \divo \bv =0}\,,
  \\
  Y := L^{p'}(\Omega)\,,\qquad  Q :=
  L^{p'}_0(\Omega),
\end{gather*}
for the most often used function spaces.

\subsection{Basic properties of the extra stress tensor} 
\label{sec:stress_tensor}
For a tensor $\bfP \in \setR^{3 \times 3} $ we denote its symmetric part by $\bP^\sym:=
\frac 12 (\bP +\bP^\top) \in \setR_\sym ^{3 \times 3}:= \set {\bfA \in \setR^{3 \times 3}
  \,|\, \bP =\bP^\top}$. The scalar product between two tensors $\bP, \bQ$ is denoted by
$\bP\cdot \bQ$, and we use the notation $\abs{\bP}^2=\bP \cdot \bP $. We assume that the
extra stress tensor $\bS$ has $(p,\para)$-structure, which will be defined now.  A
detailed discussion and full proofs of the following results can be found
in~\cite{die-ett,dr-nafsa}.
\begin{assumption}
  \label{ass:1}
  We assume that ${\bS\colon \setR^{3 \times 3} \to \setR^{3 \times 3}_\sym }$ belongs to
  $C^0(\setR^{3 \times 3};\setR^{3 \times 3}_\sym )\cap C^1(\setR^{3 \times 3}\setminus
  \{\bfzero\}; \setR^{3 \times 3}_\sym ) $, satisfies ${\bS(\bP) = \bS\big (\bP^\sym \big
    )}$, and $\bS(\mathbf 0)=\mathbf 0$. Moreover, we assume that $\bS$ has {\rm
    $(p,\para)$-structure}, i.e., there exist $p \in (1, \infty)$, $\para\in [0,\infty)$,
  and constants $C_0, C_1 >0$ such that
   \begin{subequations}
     \label{eq:ass_S}
     \begin{align}
       \sum\nolimits_{i,j,k,l=1}^3 \partial_{kl} S_{ij} (\bP) Q_{ij}Q_{kl} &\ge C_0 \big
       (\para +|\bP^\sym|\big )^{{p-2}} |\bQ^\sym |^2,\label{1.4b}
       \\
       \big |\partial_{kl} S_{ij}({\bP})\big | &\le C_1 \big (\para +|\bP^\sym|\big
       )^{{p-2}},\label{1.5b}
     \end{align}
   \end{subequations}
   are satisfied for all $\bP,\bQ \in \setR^{3\times 3} $ with $\bP^\sym \neq \bfzero$ and
   all $i,j,k,l=1,\ldots, 3$.  The constants $C_0$, $C_1$, and $p$ are called the {\em
     characteristics} of $\bfS$.
\end{assumption}
\begin{remark}
  {\rm 
    We would like to emphasize that, if not otherwise stated, the constants in the
    paper depend only on the characteristics of $\bfS$, but are independent of
    $\delta\geq0$.  
}
\end{remark}
Another important tool are {\rm shifted N-functions}\footnote{A function
$\psi\colon \setR^{\ge 0} \to \setR^{\ge 0} $ is called {\em
  N--function} if it is a continuous, convex function such that
$\lim_{t\to 0} \frac{\psi(t)}t =0$, $\lim_{t\to \infty}
\frac{\psi(t)}t =\infty$, and $\psi (t) >0 $ for $t>0$. } $\set{\phi_a}_{a \ge 0}$,
cf.~\cite{die-ett,die-kreu,dr-nafsa}. Defining for $t\geq0 $ a special N-function $\phi$
by
\begin{align} 
  \label{eq:5} 
  \varphi(t):= \int\limits _0^t \varphi'(s)\, ds\qquad\text{with}\quad
  \varphi'(t) := (\delta +t)^{p-2} t\,,
\end{align}
we can replace $C_i \big (\para +|\bP^\sym|\big )^{{p-2}}$ in the right-hand side
of~\eqref{eq:ass_S} by $\widetilde C_i\,\varphi'' \big (|\bP^\sym|\big
)$, with some constants $\widetilde C_i>0$, $i=0,1$. Next,
the shifted N-functions are defined for $t\geq0$ by
\begin{align}
  \label{eq:phi_shifted}
  \varphi_a(t):= \int\limits _0^t \varphi_a'(s)\, ds\qquad\text{with }\quad
  \phi'_a(t):=\phi'(a+t)\frac {t}{a+t}.
\end{align}
Note that $\phi_a(t) \sim (\delta+a+t)^{p-2} t^2$ and
that the complementary function satisfies $(\phi_a)^*(t) \sim
((\delta+a)^{p-1} + t)^{p'-2} t^2$.  Moreover, the {N-func\-}tions
$\phi_a$ and $(\phi_a)^* $ satisfy the
$\Delta_2$-condition\footnote{An N-function $\psi$ satisfies the
  $\Delta_2$-condition if there exists a constant $K$ such that $\psi
  (2t)\le K\, \psi (t)$ for all $t\ge 0$. The smallest such constant
  is denoted by $\Delta_2(\psi)$.} uniformly
with respect to $a\ge 0$, i.e., $\Delta_2(\phi_a) \le c\, 2^{\max
  \{2,p\}}$ and $\Delta_2((\phi_a)^*) \le c\, 2^{\max
  \{2,p'\}}$, respectively. We will use also Young's inequality: for all
$\varepsilon >0$ there exists $c_\epsilon>0 $, such that for all $s,t,a\geq 0$ it holds
 \begin{align}
   \label{ineq:young}
   \begin{split}
     ts&\leq \epsilon \, \phi_a(t)+ c_\epsilon \,(\phi_a)^*(s)\,,
     \\
     t\, \phi_a'(s) + \phi_a'(t)\, s &\le \epsilon \, \phi_a(t)+ c_\epsilon \,\phi_a(s).
   \end{split}
 \end{align}

 Closely related to the extra stress tensor $\bS$ with $(p,\delta)$-structure is the
 function $\bF\colon\setR^{3 \times 3} \to \setR^{3 \times 3}_\sym$ defined through
\begin{align}
  \label{eq:def_F}
  \bF(\bP):= \big (\para+\abs{\bP^\sym} \big )^{\frac {p-2}{2}}{\bP^\sym } \,.
\end{align}
In the following lemma we recall several useful results, which will be frequently used in
the paper. The proofs of these results and more details can be found
in~\cite{die-ett,dr-nafsa,die-kreu,bdr-phi-stokes}.
\begin{proposition}
  \label{lem:hammer}
  Let $\bfS$ satisfy Assumption~\ref{ass:1}, let $\phi$ be defined
  in~\eqref{eq:5}, and let $\bfF$ be defined in~\eqref{eq:def_F}.
  \begin{itemize}
  \item [\rm (i)] For all $\bfP, \bfQ \in \setR^{3 \times 3}$ 
    \begin{align*}
        \big({\bfS}(\bfP) - {\bfS}(\bfQ)\big) \cdot \big(\bfP-\bfQ
        \big) &\sim \bigabs{ \bfF(\bfP) - \bfF(\bfQ)}^2,
        \\
        &\sim \phi_{\abs{\bfP^\sym}}(\abs{\bfP^\sym - \bfQ^\sym}),
        \\
        &\sim \phi''\big( \abs{\bfP^\sym} + \abs{\bfQ^\sym}
        \big)\abs{\bfP^\sym - \bfQ^\sym}^2,
      \\
      \bfS(\bfQ) \cdot \bfQ &\sim \abs{\bfF(\bfQ)}^2 \sim
        \phi(\abs{\bfQ^\sym}),
      \\
      \abs{\bfS(\bfP) - \bfS(\bfQ)} &\sim
        \phi'_{\abs{\bfP^\sym}}\big(\abs{\bfP^\sym - \bfQ^\sym}\big).
    \end{align*}
  The constants depend only on the characteristics of $\bfS$.
\item [\rm (ii)] For all $\epsilon>0$, there exist a constant $c_\epsilon>0$ (depending
  only on $\epsilon>0$ and on the characteristics of $\bfS$) such that for all $\bu,
  \bv,\bw \in X$ we have
  \begin{align*}
    &\big ( {\bfS(\bD\bfu) - \bfS(\bD\bfv)},{\bD\bfw - \bD
      \bfv}\big )
      \leq \epsilon\, \norm{\bfF(\bD\bfu) - \bfF(\bD\bfv)}_2^2
      +c_\epsilon\,  \norm{\bfF(\bD\bfw) - \bfF(\bD\bfv)}_2^2\,,
    \\
      &\big ( {\bfS(\bD\bfu) - \bfS(\bD\bfv)},{\bD\bfw - \bD
      \bfv}\big )
      \leq \epsilon\, \norm{\bfF(\bD\bfw) - \bfF(\bD\bfv)}_2^2
      +c_\epsilon\,  \norm{\bfF(\bD\bfu) - \bfF(\bD\bfv)}_2^2\,,
\end{align*}
  and for all $\bfP,\bfQ\in\setR^{3 \times 3}_\sym$, $t\geq 0$ it holds
  \begin{align*}
    \phi_{\abs{\bfQ}}(t)&\leq c_\vep\, \phi_{\abs{\bfP}}(t)
    +\vep\, \abs{\bfF(\bfQ) - \bfF(\bfP)}^2,
    \\
    (\phi_{\abs{\bfQ}})^*(t)&\leq c_\vep\, (\phi_{\abs{\bfP}})^*(t)
    +\vep\, \abs{\bfF(\bfQ) - \bfF(\bfP)}^2\,.
  \end{align*}
  \item [\rm (iii)]   For all
  $\bfH \in L^p(\Omega)$ there holds  
  \begin{align*}
    \int\limits_\Omega \abs{\bfF(\bfH) - \mean{\bfF(\bfH)}_\Omega}^2\, dx \sim  \int\limits_\Omega
    \abs{\bfF(\bfH) - \bfF(\mean{\bfH}_\Omega)}^2\, dx,
  \end{align*}
  with constants depending only on~$p$.
  \end{itemize}
\end{proposition}

Let us recall the following result, taken from~\cite[Lemma
8]{dr-7-5}, \cite[Lemma 4.1]{bdr-7-5}, which is valid for $p \le 2$. 

\begin{lemma}\label{july_1}
  Let $\bS$ satisfy Assumption~\ref{ass:1}  with $p \in (1,2]$ and
  $\delta \in [0,\infty)$. Then, there exists a constant $c$, depending only on the characteristics
  of $\bS$, such that  for sufficiently smooth $\bfu$,
  $\bfv$  there holds
  \begin{align*}
    \norm{\bfF(\bfD\bfu)-\bfF(\bfD\bfv)}^2_2 \geq c
     \,\big (\delta+ \norm{\bfD \bfu}_p +    \norm{\bfD \bfu - \bfD
    \bfv}_p\big )^{p -2} \, \norm{ \bfD \bfu - \bfD \bfv}^2_p\,.
  \end{align*}

\end{lemma}

\subsection{Discretizations}\label{sec:disc}
For the time-discretization, given $T>0$ and $M\in \setN$, we define the time step size as
$\kappa:=T/M>0$, with the corresponding net $I^M:=\{t_m\}_{m=0}^M$, $t_m:=m\,\kappa$. We
use the notation $I_m:=(t_{m-1},t_m]$, with $m=1,\ldots,M$.  For a given sequence
$\{\bv^m\}_{m=0}^M$ we define the backward differences quotient as
\begin{equation*}
  d_t \bv^m:=\frac{\bv^m-\bv^{m-1}}{\kappa}\,, \quad m=1,\ldots,M.
\end{equation*}


The proof of the main result uses the following
modification of Gronwall's lemma.

\begin{lemma}
  \label{lem:discrete_Gronwall_lemma}
  Let $1<p\leq 2$ and $T\in (0,\infty)$.  For $M\in \setN$ and $h>0$ let be given non-negative sequences
  $\{a_m(h)\}_{m=0}^M$, $\{b_m(h)\}_{m=0}^M $, $\{r_m(h,\kappa)\}_{m=1}^M$,
  $\{s_m(h,\kappa)\}_{m=1}^M$, $\{\rho_m(h,\kappa)\}_{m=1}^M$ and
  $\{\sigma_m(h,\kappa)\}_{m=1}^M$, where $\kappa:=\frac TM$. Assume
  that there exists  ${\mu_0, \widehat \kappa>0}$ such that for all
  $0<h<1/\sqrt{\mu_0}$ and all $0< \kappa < \widehat \kappa$ there holds:
  \begin{equation}
    \label{eq:r_m}
    \begin{gathered}
      \big(a_0 (h)\big)^{2}\leq \mu_0\,h^{2}, \qquad
      \big(b_0 (h)\big)^{2} \leq \mu_0\,h^{2},
      \\
      \ksumo \big(r_m (h,\kappa)\big)^{2} \leq
      \mu_0\,h^2,\qquad\ksumo \big(s_m (h,\kappa)\big)^{2} \leq \mu_0\,h^2,
      \\
      \ksumo \big(\rho_m (h,\kappa)\big)^{2}\leq
      \mu_0\,\kappa^2,\qquad\ksumo \big(\sigma_m (h,\kappa)\big)^{2}\leq \mu_0\,\kappa^2.
    \end{gathered}
  \end{equation}
  Further, let there exist constants $\mu_1,\mu_2,\mu_3>0$,
  $ \Lambda>0$, and some $0<\theta\leq 1$ such that for some $\lambda
  \in [0,\Lambda]$ the following two inequalities are satisfied for
  all $0<h<1/\sqrt{\mu_0}$,  all $0< \kappa < \widehat \kappa$ and all
  $m=1,\ldots, M$:\footnote{Here we use the convention that for
    $\lambda=b_m=0$ we set $(\lambda+b_m)^{p-2}b_m^2=0$.}
  \begin{align}
    \label{discrete_Gronwall_bis}
    &d_t a_m^2+\mu_1(\lambda+b_m)^{p-2}b_m^2\leq b_m r_m+
    b_m \rho_m+\mu_2 b_{m-1}b_m+s_m^2 +\sigma_m^2,
    \\
    \label{discrete_Gronwall_ter}
    &d_t a_m^2+\mu_1(\lambda+b_m)^{p-2}b_m^2\leq b_m r_m+b_m \rho_m+\mu_3
    b_{m}b_{m-1}^{1-\theta}a_{m-1}^\theta+s_m^2 +\sigma_m^2.  
  \end{align}
  %
  %
  Then, there exist constants $\overline{\mu_0}, \overline{\kappa} >0$, 
  and $\mu_4,\,\mu_5>0$, independent of $\lambda$, such that for all
  $\kappa, h>0$ satisfying $\kappa<\overline{\kappa}$  and 
  $h^{2}<\overline{\mu_0}\,\kappa$ there holds 
  \begin{align}
    &\max_{0\leq m\leq M} b_m\leq 1,
    \\
    \label{Gronwall_discrete2}
    &\max_{0\leq m\leq M} a_m^2+{\mu_1(1+\Lambda)^{p-2}}
    \ksum b_m^2\leq \mu_4\,\big (h^{2} +\kappa^2\big )\,\textrm{exp}(2 \mu_{5}\kappa\,M).
  \end{align}
\end{lemma}

\begin{proof}
  This result is a small modification of the corresponding results
  in~\cite{bdr-7-5-num,bdr-3-2}, and can be proved in the same way. 
\end{proof}

The following result will
be used frequently in the sequel. 
\begin{lemma}
  \label{lem:Bochner-lemma}
 Assume that 
  \begin{equation*}
    f, {\partial_t f}\in L^2(I;X)\,,
  \end{equation*}
where $(X,\|\cdot\|_{X})$ is a Banach space. Then, for all $\tau _m \in I_m$, $m=1,\ldots, M$, it holds 
\begin{equation}
  \label{eq:estimate-time-derivative}
  \begin{aligned}
    \ksumo\;\Int\|f(s)-f(\tau_{m})\|_X^2\,ds &\leq
    \kappa^2\,\bignorm{{\partial_t f}}_{L^2(I;X)}^2\,.
  \end{aligned}
\end{equation}
\end{lemma}
\begin{proof}
  The assertion is proved in \cite[Lemma 3.1]{br-parabolic} in the
  special case $\tau _m =t_m$, $m=1,\ldots, M$. The general
  case follows exactly in the same way.
\end{proof}

For the spatial discretization we denote by 
$\mathcal{T}_h$ a family of shape-regular triangulations, consisting of
$3$-dimensional \farb{closed} simplices $K$. We denote by $h_K$ the diameter of $K$ and by $\rho_K$ the
supremum of the diameters of inscribed balls. We assume that $\mathcal{T}_{h}$ is
non-degenerate, i.e., there exists a constant $\gamma_0>0$ such that $\max _{K \in \mathcal{T}_{h}} \frac {h_K}{\rho_K}\le \gamma_0$.  The
global mesh-size $h$ is defined by $h:=\max _{K \in \mathcal{T}_h}h_K$.  Let $S_K$ denote the
neighborhood of~$K$, i.e., $S_K$ is the union of all simplices of~$\mathcal{T}_{h}$
intersecting~$K$.  By the assumptions we obtain that $|S_K|\sim |K|$ and that the number of
patches $S_K$ to which a simplex belongs are bounded uniformly in both $h>0$ and $K\in \mathcal{T}_{h}$.

We denote by ${\mathcal P}_k(\mathcal{T}_{h})$, with $k \in \setN_0:=\setN\cup\{0\}$, the
space of scalar or vector-valued 
functions, which are polynomials of degree at most $k$ on each $K\in \mathcal{T}_{h}$.
Given a triangulation $\mathcal{T}_h$ of $\Omega$ with the above properties and given
$r_0, r_{1}, s_0\in\setN_{0}$, with $r_0 \le r_1$, we define
\begin{align*}
  X_{h}:=\left\{\bv_h\in X\fdg \bv_h\in \mathcal{P}\right\}
  \qquad \text{and}\qquad
  Y_{h}:=\left\{\eta_h\in Y\fdg \eta_h\in \mathcal{P}_{s_0}(\mathcal{T}_h)\right\},
\end{align*}
with ${\mathcal
  P}_{r_0}(\mathcal{T}_h) \subset \mathcal{P} \subseteq {\mathcal P}_{r_1}(\mathcal{T}_h)$.
Note that there exists a constant $c=c(r_1, \gamma_0)$
such that for all $\bv_h \in X_h$, $K \in \mathcal{T}_h$, $j\in
\setN_0$, and all $x \in K$ holds
\begin{align}
  \label{eq:inverse}
  \abs{\nabla ^j \bv_h(x)}\le c\, \dashint_K \abs{\nabla ^j \bv_h(y)}\, dy\,.
\end{align}
For the weak formulation of discrete problems we
use the following function spaces
\begin{equation*}
    V_{h}:=V\cap X_{h}\,, \qquad Q_h:=Q\cap Y_h\,.
\end{equation*}

We also need some numerical interpolation operators.  Rather than
working with a specific interpolation operator we make the following
assumptions:
\begin{assumption}
  \label{ass:proj}
   We assume that $r_0=1$ 
  and that there exists a linear projection operator
  $\Pidiv \colon X \to X_h$ which 
  \begin{enumerate}
  \item is locally $W^{1,1}$-stable, i.e., for all $\bfw
    \in X$ and $K \in \mathcal{T}_h$ there holds
    \begin{align}
      \label{eq:Pidivcont}
      \dashint_K \abs{\Pidiv \bfw}\,dx &\leq c \dashint_{S_K}\!
      \abs{\bfw}\,dx + c \dashint_{S_K}\! h_K \abs{\nabla \bfw}\,dx
     \,;
    \end{align}
  \item preserves zero boundary values, i.e., $\Pidiv (V) \subset V_h$;
    \item preserves divergence in the $Y_h^*$-sense, i.e.,  for all
      $\bfw \in X$ and  $\eta_h \in
      Y_h$ there holds 
    \begin{align*}
      \hskp{\divergence \bfw}{\eta_h} &= \hskp{\divergence \Pidiv
        \bfw}{\eta_h}\,.
    \end{align*}
  \end{enumerate}
\end{assumption}

\begin{assumption}
  \label{ass:PiY}
  We assume that $Y_h$ contains the constant functions, i.e. that
  $\setR \subset Y_h$, and that there exists a linear projection
  operator $\PiY\,:\, Y \to Y_h$ which is locally $L^1$-stable, i.e.,
  for all $q \in Y$ and $ K \in \mathcal{T}_h$ there holds
  \begin{align*}
    \dashint_K \abs{\PiY q}\,dx &\leq c\, \dashint_{S_K} \abs{q}\,dx
   \,.
  \end{align*}
\end{assumption}%

The existence of a projection operator $\Pidiv$ as in
Assumption~\ref{ass:proj} is known (among others) for the Taylor-Hood,
the Crouzeix-Raviart, and the MINI element in dimensions two and
three; the Cl\'ement interpolation operator satisfies
Assumption~\ref{ass:PiY}.  For a discussion and consequences of these
assumptions we refer to~\cite[Sec.~3.2]{bdr-phi-stokes},
\cite[Appendix]{br-parabolic}, and \cite[Sec.~4,5]{dr-interpol}. We
collect in the next two propositions the properties of the projection
operators, which are relevant for the present paper.
\begin{proposition}\label{prop:Ph}
  Let $\Pidiv $ satisfy Assumption~\ref{ass:proj}. 
  \begin{enumerate}
  \item [\rm (i)]   Let $\bF(\bD \bv) \in W^{1,2}(\Omega)$. Then, there exists a constant
    $c=c(p,r_1,\gamma_0)$ such that 
  \begin{align*}
    \norm{\bF(\bD\bv) -\bF(\bD \Pidiv \bv)}_2
    &\le c\, h\, \norm{\nabla \bF(\bD \bv ) }_2. 
  \end{align*}
  \item [\rm (ii)] Let $r\in(1,\infty)$. Then, there exists a constant
    $c=c(r,r_1,\gamma_0)$ such that for all $\bv \in
    W^{1,r}(\Omega)$ and all $\bw \in
    W^{2,r}(\Omega)$  holds
    \begin{align*}
      \|\bv-\Pidiv \bv\|_r+h\,\|\nabla \Pidiv \bv \|_r&\le c\,  h\,\|\nabla \bv\|_r\,,
      \\
      \|\nabla \bw- \nabla \Pidiv \bw\|_r &\le c\,  h\,\|\nabla ^2\bw\|_r\,.
\end{align*}
  \item [\rm (iii)] Let $r\in[1,2]$ and let $\ell =1$ or $\ell =2$ be such
    that $W^{\ell,r}(\Omega) \vnor \vnor L^2(\Omega)$. Then, there exists a constant
    $c=c(r,\ell,r_1,\gamma_0)$ such that for all $\bv \in
    W^{\ell,r}(\Omega)$ holds
    \begin{align*}
      \|\bv-\Pidiv \bv\|_2\le c\,  h^{\ell +3( \frac 12 -\frac 1r)}\,\|\nabla^\ell \bv\|_r\,.
    \end{align*}
  \item [\rm (iv)] Let $\bF(\bD \bv) \in
    W^{1,2}(\Omega)$ and  $\bF(\bD \bw) \in
    L^{2}(\Omega)$. Then, there exists a constant
    $c=c(p,r_1,\gamma_0)$ such that 

\begin{align*}
  \begin{aligned}
    &\int\limits_\Omega \phi_{|\bD\bfv|}\big(\big|\bD \Pidiv {\bv}
    -\bD \Pidiv \bfw \big|\big)\,\dx
\leq c \,h^2\|\nabla\bF(\bD\bv)\|^2_2
    +c\,\|\bF(\bD\bv)-\bF(\bD\bw)\|^2_2 \,.
      \end{aligned}
    \end{align*}
\end{enumerate}
\end{proposition}
\begin{proof}
  The first two assertions are proved,
  e.g.,~in~\cite{bdr-phi-stokes,dr-interpol}. The last two assertions 
  are proved in \cite{br-parabolic}. 
\end{proof}
\begin{proposition}\label{pro:PY}
  Let $\PiY$ satisfy Assumption \ref{ass:PiY}. Let $\psi$ be an
  N-function satisfying the $\Delta_2$-condition. Then, there exists a
  constant $c=c(\gamma_0, \Delta_2(\psi))$ such that for all
  sufficiently smooth functions and all $K \in \mathcal T_h$ there
  holds
  \begin{align*}
    \int\limits_K \psi(\abs{\PiY q})\, dx
    &\le c     \int\limits_{S_K} \psi(\abs{q})\, dx\,,
    \\
    \int\limits_K \psi(\abs{q-\PiY q})\, dx
    &\le c     \int\limits_{S_K} \psi(h_K\abs{\nabla q})\, dx\,.
  \end{align*}
\end{proposition}
\begin{proof}
  The assertions are proved in~\cite{dr-interpol}.
\end{proof}
\subsection{Main results}\label{sec:main}
Before we formulate the main result, proving optimal convergence rates for the error
between the solution $\bu$ of the continuous problem~\eqref{eq:pfluid} and the discrete
solution $\{\bfuMh \}_{m=0}^{ M}$ of the space-time
discretization~\eqref{eq:fdiscr-var}, we discuss the existence of
solutions of \eqref{eq:pfluid} and \eqref{eq:fdiscr-var}.

The existence of global weak solutions
$\bfu \in L^\infty(I;L^2(\Omega))\cap L^p(I;V_{\divo} )$ of
\eqref{eq:pfluid} for large data is proved in \cite{die-ru-wolf} for
$\delta\ge 0$ and $p>\frac{6}{5}$. The existence of a locally in time,
unique strong solution
$\bfu \in L^r(I';W^{2,r}(\Omega)) \cap L^p(I';V_{\divo} ) \cap
W^{1,r}(I';L^r(\Omega))$,
$q \in L^r(I';W^{1,r}(\Omega)) \cap L^{r}(I';L^r_0(\Omega))$ for any
$r\in (5,\infty)$ and some $I':=(0,T')$ with $ {T'\in (0,T)}$ of
\eqref{eq:pfluid} for large data is proved in \cite{pruess} for
$\delta > 0$ and $p>1$. 

The existence of a unique discrete solution
$\{\bfumh\}_{m=1}^M, \{\piMh\}_{m=1}^M$ can be inferred in the
following way.  Setting
$V_h(0):=\{\bfw_h\in V_h \fdg \hskp{\divergence\bfw_h}{\eta_h}=0 \
\forall\,\eta_h\in Y_h\}$, one uses Brouwer's fixed point theorem and
the properties of the extra stress tensor $\bS$ and the convective term to
show the existence of $\bfuMh \in V_h(0)$ satisfying
\eqref{eq:fdiscr-var}$_1$ for all $\bfv_h \in V_h(0)$ (cf.~\cite[Lemma
7.1]{pr}). This solution is unique due to the semi-implicit
discretization of the convective term and the monotonicity of
$\bS$. Moreover,  testing
\eqref{eq:fdiscr-var} with $\{\bfu^m_h \}_{m=1}^{ M}$ yields the energy estimate
\begin{equation*}
  \max_{m=1,\dots,M}\|\bfu^m_h\|^2_2+\ksumo\|\bF(\bD\bfu^m_h)\|^2_2\leq
  C(\bu_0, \ff)\,.
\end{equation*}
 The properties of  the projection operators $\Pidiv$ and
$\PiY$ ensure the validity of the discrete inf-sup condition
(cf.~\cite[Lemma~4.1]{bdr-phi-stokes}). Since the 
divergence is a closed, surjective, linear operator from $V_h$ onto
$Q_h$ and since the gradient is the annihilator of
the kernel of the divergence, these  ensure the existence of a unique $\piMh
\in Q_h$ such that $\{\bfuMh\}_{m=1}^{ M}, \{\piMh\}_{m=1}^{ M}$ satisfy \eqref{eq:fdiscr-var}. 

Now we can formulate our main result.

\begin{theorem}
\label{thm:theorem-parabolic}
Let the extra stress tensor $\bS$ in~\eqref{eq:pfluid} have
$(p,\delta)$-structure for some $p\in(\frac 65,2]$, and some
$\delta\in[0,\infty)$ fixed but arbitrary.  Let $\Omega\subset\setR^3$
be a bounded polyhedral domain with Lipschitz continuous boundary, and
$I=(0,T)$, $T\in (0,\infty)$, be a finite time interval. Assume that
$\ff \in 
{W^{1,2}(I;L^2(\Omega))}$,
$\bu_0 \in W^{2,p}(\Omega)\cap V_{\divo}$ and that the solution
$(\bu, q)$ of~\eqref{eq:pfluid} satisfies~\eqref{eq:cont-var} and
\begin{equation}\label{eq:reg-ass} 
  {\bF(\bD \bu)}\in W^{1,2}(I \times\Omega) \qquad\text{and}\qquad  q \in
  L^{p'}(I;L^{p'}_0(\Omega)\cap W^{1,p'}(\Omega))\,.
\end{equation}
Let the space $V_h$, $h>0$, be defined as above with $r_0=1$ and let
$\{\bfuMh\}_{m=1}^{ M}$,  $\{\piMh\}_{m=1}^{ M}$  be the unique
solution of~\eqref{eq:fdiscr-var}.
Then, for $p\in(\frac 85,2]$
there exists $\kappa_0\in (0,1]$ such that, for given $h\in ( 0, 1)$ and 
$\kappa \in ( 0, \kappa_0)$ satisfying
\begin{equation}
  \label{eq:compatbility}
  h^{4/p'}\leq \sigma_0\, \kappa,
\end{equation}
for some $\sigma_{0}>0$, the following error estimate holds true 
\begin{equation}\label{eq:err}
  \max_{m=1,\ldots,M}
  \|\bfuMh-\bfu(t_{m})\|_2^2+\ksumo\|\bF(\bD\bfuMh)-\bF(\bD\bfu(t_{m}))\|_2^2\leq
  c\,(h^2+\kappa^2),
\end{equation}
with a constant $c$ depending only on the characteristics of $\bS$,
$\norm{\bF(\bD \bu)} _{W^{1,2}(I \times\Omega)}$,
$\norm{\partial _t\ff}_{L^2(I;L^2(\Omega)}$, $\norm{\bu_0}_{2,p}$,
$\norm{\nabla q}_{L^{p'}(I\times \Omega)}$,
$\gamma_0$, $r_1$, $\delta$, and $\sigma_0$.
\end{theorem}

\subsection{Comparison with previous results and observation on the requested regularity.}
Here we compare the new result with previous ones on the discretization
of generalized non-Newtonian fluids (and general parabolic equations and systems). We also
discuss briefly the regularity we are assuming on the continuous solution.

\smallskip

Concerning previously proved error-estimates we can mainly recall the following facts:

  (i) Problem \eqref{eq:pfluid}, in the case of space periodic
  boundary conditions, has been treated
  in~\cite{bdr-3-2,dpr1,dpr2,pr}. In \cite{bdr-3-2} the same optimal error estimate as in Theorem
  \ref{thm:theorem-parabolic} is proved under slightly stronger
  assumptions on the regularity of the solution $\bu$ of
  \eqref{eq:pfluid}, for $p \in (\frac32,2]$. 

Problem \eqref{eq:pfluid}
  without the convective term $[\nabla \bu]\bu$,  in the case of
  homogeneous Dirichlet boundary conditions, has been treated
  in~\cite{sarah-phd} for $p \in (\frac 65,\infty)$. It is shown there
  that \eqref{eq:err} holds with the right-hand side replaced by
  $c \,(h^{\min \{2,\frac 4p\}} +\kappa^2)$.  The proofs of  these
  results are based on intermediate semi-discrete problems, for which
  a certain regularity has to be proved, to obtain the desired optimal
  convergence rates. This in fact limits the results
  in~\cite{bdr-3-2,dpr1,dpr2,pr} to the case of space periodic
  boundary conditions. Here, we avoid such (technical) problems by proving the
  error estimate directly without using intermediate semi-discrete
  problems. The same perspective is also taken in the recent papers
\cite{br-parabolic,breit-mensah}. In \cite{breit-mensah} problem
  \eqref{eq:pfluid} is treated for a nonlinear operator $\bS$
  depending on the full gradient $\nabla \bu$, having
  $(p(\cdot,\cdot),\delta)$-structure, with a variable exponent
  $p(\cdot,\cdot)$, but without convective term and without the
  solenoidality condition. In this situation the error estimate
  \eqref{eq:err} with the right-hand side replaced by
  $c\, (h^{2\alpha_x} +\kappa^{2\alpha_t})$ is proved if the variable
  exponent belongs to
  $C^{\alpha_x,\alpha_t}(\overline{I\times \Omega})$ for appropriate
  $\alpha_x, \alpha_t \in (0,1]$ and an additional CFL-condition
  $\kappa ^r \le c\, \inf_{K \in \mathcal T_h}h_k$, with some
  $r\ge \frac {1+2\alpha_t}{2\alpha_x}$, is satisfied. In
  \cite{br-parabolic} problem \eqref{eq:pfluid} is treated without
  convective term and without the solenoidality condition. There the
  error estimate \eqref{eq:err} is
  proved under the same conditions as in Theorem~\ref{thm:theorem-parabolic}. \\[-3mm]

  (ii) In the recent paper \cite{breit-lars-etal} a different approach
  is used to treat the unsteady $p$-Laplace problem, i.e., problem
  \eqref{eq:pfluid} without convective term and without the
  solenoidality constraint.  By using the $L^{2}$-projection operator
  instead of the Scott--Zhang operator (cf.~\cite{zhang-scott}),
  satisfying Assumption \ref{ass:proj}, the error estimate
\begin{equation}\label{eq:err-lars-dominic}
  \begin{aligned}
   & \max_{m=1,\ldots,M}
    \|\bfuMh-\dashint\nolimits_{t_{m-1}}^{t_{m+1}}\bfu(s)\,ds\|_2^2
    \\
    & \qquad +\ksumo\int_{t_{m-1}}^{t_{m+1}}\|\bF(\bD\bfuMh)-\bF(\bD\bfu(s))\|_2^2\,ds
    \leq c\,(h^{2\alpha_{s}}+\kappa^{2\alpha_{t}}),
  \end{aligned}
\end{equation}
is proved in \cite{breit-lars-etal} without a coupling condition
between $h$ and $\kappa$. The removal of the coupling is due to the
usage of the $L^{2}$-projection which in the treatment of the time
derivative does not produce terms needing a coupling
(cf.~estimate~\eqref{eq:time}). However, the treatment of the
nonlinear operator $\bS$ is subtle and requires delicate estimates,
which result in the different error
estimate~\eqref{eq:err-lars-dominic} compared to the error estimate \eqref{eq:err}. The estimate
\eqref{eq:err-lars-dominic} is proved under the assumption that
\begin{equation*}
  \begin{aligned}
    {\bF(\nabla\bu)}& \in L^{2}(I;N^{\alpha_{x},2}(\Omega))\cap
    N^{\alpha_{t},2}(I;L^{2}(\Omega))
    \\
    \bu&\in L^{2}(I;N^{\alpha_{x},2}(\Omega))
  \end{aligned}
\end{equation*}
for $\frac{1}{2}<\alpha_{t}\leq1$ and $0<\alpha_{x}\leq1$, where
$N^{\alpha,2}$ are appropriate Nikol'ski\u{\i} spaces with
differentiability $\alpha\in(0,1]$ 
(cf.~\cite{breit-lars-etal}). Even for $\alpha_{x}=\alpha_{t}=1$
estimate \eqref{eq:err-lars-dominic} differs from~\eqref{eq:err} since
it contains time averages of the error rather than a point-wise error.
The usage of limited
regularity in the time-variable and time averaging of the error is motivated by a similar analysis performed
for stochastic parabolic equations in~\cite{bhl}.\\[-3mm]

  (iii) In \cite{sueli-tscherpel,tscherpel-phd} the convergence of a
  fully implicit space-time discretization (without convergence rate
  but also with no assumptions of smoothness of the limiting problem)
  of the problem \eqref{eq:pfluid} in the case of homogeneous
  Dirichlet boundary conditions is proved for $p>\frac{11}5$ and even
  for $p>\frac 65$ for a regularized problem. The
  convergence of the same numerical scheme~\eqref{eq:fdiscr-var}
  towards a weak solution has been recently proved
  in~\cite{alex-rose-nonconform} for $p>\frac{11}5$. In fact, in
  \cite{alex-rose-nonconform} the convergence of a general quasi
  non-conforming Rothe--Galerkin scheme in the context of evolution
  problems with Bochner pseudo-monotone operators is proved
  (cf.~\cite{br-fully} for the treatment of a conforming Rothe--Galerkin
  scheme in the context of evolution problems with Bochner
  pseudo-monotone operators). 

\medskip

Let us now discuss the \textit{natural} regularity assumption \eqref{eq:reg-ass}. The
assumption \eqref{eq:reg-ass}$_1$ is natural in the sense that it is
the one obtained from the extra stress tensor $\bS$ if formally tested
with $-\Delta \bu$ and $\partial _t^2 \bu$. The existence of solutions
satisfying \eqref{eq:reg-ass}$_1$ is proved rigorously for problem
\eqref{eq:pfluid} in the case of periodic boundary conditions locally
in time in \cite{bdr-7-5, dr-7-5},  for $p>\frac 75$ and large
data. The situation for Dirichlet boundary conditions is more
complicated. The existence of a locally in time unique strong solution
${\bfu \in L^r(I';W^{2,r}(\Omega)) \cap L^p(I';V_{\divo} ) \cap
  W^{1,r}(I';L^r(\Omega))}$,
$q \in L^r(I';W^{1,r}(\Omega)) \cap L^{r}(I';L^r_0(\Omega))$ for any
$r\in (5,\infty)$ and some $I':=(0,T')$ with $ T'\in (0,T)$ is proved
in \cite{pruess} for large data. This regularity implies, using
parabolic embedding theory (cf.~\cite[Appendix]{dr-7-5}), that
$\bF(\bD \bu) \in L^2(I';W^{1,2}(\Omega))$ and that
$\bu, \nabla \bu \in C(\overline {I'\times \Omega})$. However, in
\cite{pruess} it is
not proved that this solution also satisfies $\partial _t\bF(\bD \bu)
\in L^2(I';L^{2}(\Omega))$. Nevertheless, one can show that the
solution from \cite{pruess} also satisfies $\partial _t\bF(\bD \bu)
\in L^2(I';L^{2}(\Omega))$, using the following auxiliary result.

\begin{proposition}\label{prop:reg}
  Let the extra stress tensor $\bS$ in~\eqref{eq:pfluid} have
  $(p,\delta)$-structure for some $p\in(\frac 65,2]$, and some
  $\delta\in[0,\infty)$.  Let
  $\Omega\subset\setR^d$, $d\ge 2$, be a bounded domain with Lipschitz
  continuous boundary, and $I=(0,T)$, ${T\in (0,\infty)}$, be a finite
  time interval. Assume that $ {\bfu_0} \in V_{\divo}$
  satisfies $\divo\bS(\bD\bu_0)\in L^2(\Omega)$, $ {\bff}\in {L^{p'}(I;L^{p'}(\Omega))} \cap
    {W^{1,2}(I;L^2(\Omega))}$  and 
  ${\bG \in C(\overline {I}; W^{1,2}(\Omega))\cap
    W^{1,p'}(I;L^{p'}(\Omega))}$.  Then, there exists a unique
  weak solution $\bv \in L^\infty(I;L^2(\Omega))\cap L^p(I;V_{\divo})$ of 
  \begin{equation}
    \label{eq:pfluid1}
    \begin{aligned}
      \partial_t\bfv-\divergence \bfS(\bfD\bfv)+\nabla q&= \bff +\divo
      \bG \qquad&&\text{in } I\times \Omega,
      \\
      \divergence\bfv&=0\qquad&&\text{in } I\times \Omega,
      \\
      \bfv(0)&=\bfu_0\qquad&&\text{in } \Omega,
    \end{aligned}
  \end{equation}
  satisfying additionally $\partial _t\bv \in L^\infty(I;L^2(\Omega))$ and
  $\partial_t \bF(\bD\bv) \in L^2(I;L^2(\Omega))$.
\end{proposition}
\begin{proof}
  This result is proved using ideas from \cite{bdr-7-5, dr-7-5}.
 Using the Galerkin method and the theory of monotone
  operators one shows that there exists a unique weak solutions
  $\bv \in L^\infty(I;L^2(\Omega))\cap L^p(I;V_{\divo})$ of
  \eqref{eq:pfluid1}. Moreover, the regularity of the data allows us to
  show 
  that we can take the time derivative of
  the Galerkin equations and test with the time derivative of the
  Galerkin solution. Straightforward manipulations show that this
  produces  an estimate, showing after a limiting process in the Galerkin parameter, the
  additional regularity stated above. 
\end{proof}

Next, by using Proposition~\ref{prop:reg} with
$\bG =\bu \otimes \bu$ (where $\bu$ is the solution from
\cite{pruess}), the monotonicity of $\bS$, and the above regularity for
$\bu$ imply that the solution from \cite{pruess} 
satisfies also  $\partial_t \bF(\bD\bu) \in
L^2(I';L^2(\Omega))$. Consequently, the unique solution $(\bu,q)$
from \cite{pruess} satisfies \eqref{eq:reg-ass} with $I $ replaced by
$I'$. \\[-3mm]

It is useful to formulate the consequences of the 
assumption $\bF(\bD\bu)\in W^{1,2}(I\times \Omega)$ in 
terms of Bochner--Sobolev spaces. From~\cite[Thm.~33]{dr-7-5} and
standard embedding results it follows that 
\begin{equation*}\label{eq:embed}
   \bF(\bD\bu) \in  C(\overline{I};L^3(\Omega)).
\end{equation*}
Using $|\bD \bfu|^{p/2}+\delta^{\frac p2}
\sim|\bF(\bD\bu)|+\delta^{\frac p2} $ and the continuity of $\bP
\mapsto \bF^{-1}(\bP)$ (cf.~\cite[Lemma 3.23]{bdr-7-5}) we get
\begin{equation} \label{eq:I-interpolation-u}
 \bfu\in C(\overline {I};W^{1,3p/2}(\Omega)).
\end{equation}
In~\cite[Lemma~4.5]{bdr-7-5} it is shown that for $p\leq2$ 
\begin{align*}
    \|\nabla^2\bfu\|_{\frac{6p}{4+p}}^2
    &\leq c\, \|\nabla \bF(\bD\bu)\|_2^2(\delta+\|\nabla\bu\|_{3p/2})^{2-p},
    \\
    \bignorm{ {\partial_t \nabla \bfu}}_{\frac{6p}{4+p}}^2
    &\leq c\, \|\partial _t \bF(\bD\bu)\|_2^2(\delta+\|\nabla\bu\|_{3p/2})^{2-p}\,.
\end{align*}
Thus, we also get 
\begin{align}
  \label{eq:regularity}
  \begin{aligned}
    \bfu &\in L^2(I;W^{2,\frac{6p}{4+p}}(\Omega)),
    \\
    {\partial_t \bfu} &\in L^2(I;W^{1,\frac{6p}{4+p}}(\Omega)),
  \end{aligned}
\end{align}
where the bounds depend only on $\|\bF(\bD\bu)\|_{W^{1,2}(I\times
  \Omega)}$ and $\delta_0$. \\[-3mm]

Finally we would like to comment on the restriction $p \in (\frac
85,2] $ in Theorem \ref{thm:theorem-parabolic} compared to the
restriction $p \in (\frac 32,2]$ in \cite[Theorem 2.6]{bdr-3-2}. Based
on the results from \cite{bdr-7-5} it is assumed in
\cite{bdr-3-2} that, additionally to \eqref{eq:reg-ass}$_1$,
the solution satisfies among other properties $\bF(\bD \bu) \in
L^{2\frac{5p-6}{2-p}}(I;W^{1,2}(\Omega))$. This results in
\begin{align}\label{ass:reg2}
  \bfu\in C(\overline {I};W^{1,r}(\Omega))\qquad  \text{ for any } r\in [1,6(p-1))\,.
\end{align}
If we also assume that \eqref{ass:reg2} holds, we can improve Theorem
\ref{thm:theorem-parabolic} to the range $p \in (\frac 32,2]$. More
precisely, we have: 
\begin{cor}\label{cor:main}
  Assume that in the situation of Theorem
  \ref{thm:theorem-parabolic} the solution $\bu$ of \eqref{eq:pfluid}
  additionally satisfies \eqref{ass:reg2}. Then, the error estimate \eqref{eq:err}
  holds for ${p \in (\frac 32,2]}$ with a constant additionally depending
  on $\norm{\bu}_{C(\overline {I};W^{1,r}(\Omega))}$, for some
  suitable  $r \in [1,6(p-1))$. 
\end{cor}

\section{Proof of the main result}
In this section we prove the error estimates from
Theorem~\ref{thm:theorem-parabolic}. To this end we need to derive an 
equation for the error and to use the discrete Gronwall lemma~\ref{lem:discrete_Gronwall_lemma} together
with approximation properties due to the regularity of the solution
and the properties of the extra stress tensor $\bS$.

In the error equation we want to use the test function $\bfuMh-\Pidiv
\bfu(t_{m})$, which belongs to the space $V_h(0)$. Thus, it is enough
to consider test functions $\bv_h$ from  $ V_h(0)$ in \eqref{eq:fdiscr-var}.
For such test functions we can replace the discrete pressure
$q_h^m$ by an arbitrary function from $Q_h$. Thus,  it follows from \eqref{eq:fdiscr-var}$_1$
that for all ${\bv _h \in V_h(0)}$, $\mu_h\in Q_h$ and $m=1,\ldots, M$ there holds
\begin{align}  
  \label{eq:mod}
    (d_t \bfuMh,\bfvM)+\hskp{\bS(\bD\bfuMh)}{\bfv_h}
    +b(\bfuMsh,\bfuMh,\bfvM)-\hskp{\mu_h }{\divo
      \bfvM}&=(\bff(t_m),\bfvM)\,.
\end{align}
We can choose in~\eqref{eq:mod}, for each $m=1,\ldots,M$,
\begin{equation*}
  \mu_h=\mu_{h}^{m}:=\Int \PiY q(t)\,dt,
\end{equation*}
%
and since $\hskp{\Int\PiY q(t) \,dt}{\divo
      \bfvM}=\Int\hskp{\PiY q(t) }{\divo
      \bfvM}\,dt$, we get that for all ${\bv _h
  \in V_h(0)}$ and   $m=1,\ldots, M$,   there holds 
\begin{align*}  
    (d_t \bfuMh,\bfvM)+\hskp{\bS(\bD\bfuMh)}{\bfv_h}
    +b(\bfuMsh,\bfuMh,\bfvM)-\Int\hskp{\PiY q(t) }{\divo
      \bfvM}\,dt&=(\bff(t_m),\bfvM),
\end{align*}
which we can re-write also as follows:
\begin{align*}  
    &(d_t \bfuMh,\bfvM)+\Int\hskp{\bS(\bD\bfuMh)}{\bfv_h}\,dt
    +\Int b(\bfuMsh,\bfuMh,\bfvM)\,dt-\Int\hskp{\PiY q(t) }{\divo
      \bfvM}\,dt
\\
&=\Int(\bff(t_m),\bfvM)\,dt,
\end{align*}
since we are averaging locally constant terms.
We subtract from this equation the retarded averages over $I_m$ of equation
\eqref{eq:cont-var} and obtain the equation for the error
\begin{align}  \label{eq:error-var}
  \begin{aligned}
    &\bighskp {d_t\big (\bfuMh-\bfu(t_{m})\big )}{\bfvM}+\Int\bighskp{
      \bS(\bD\bfuMh)-\bS(\bD\bfu(t))}{\bD\bfvM} \,dt
    \\
    &\quad +\Int b(\bfuMsh,\bfuMh,\bfvM) - b(\bu(t),\bu(t),\bfvM)\, dt
    -\Int \hskp{\PiY q(t) -q(t)}{\divo \bfvM}\,dt
    \\
    &=\Int (\ff (t_m)-\bff(t),\bfvM )\, dt,
  \end{aligned}
\end{align}
valid for all $ m=1,\ldots, M$ and all $\bfvM\in V_h(0)$.

Choosing now the legitimate test function $\bfvM=\bfuMh -\Pidiv
\bu(t_m) \in V_h(0)$ we finally get for all
$m=1,\ldots, M$
\begin{align}  \label{eq:err-fin}
  \begin{aligned}
    &\bighskp {d_t\big (\bfuMh-\bfu(t_{m})\big )}{\bfuMh -\Pidiv
      \bu(t_m) }
    \\
    &\quad +\Int\bighskp{ \bS(\bD\bfuMh)-\bS(\bD\bfu(t))}{\bD \bfuMh -
      \bD\Pidiv \bu(t_m) } \,dt
    \\
    &\quad +\Int b(\bfuMsh,\bfuMh, \bfuMh -\Pidiv \bu(t_m)) - b(\bu(t),\bu(t), \bfuMh
    -\Pidiv \bu(t_m) )\, dt
    \\
    &\quad -\Int \bighskp{\PiY q (t) - q(t)}{\divo
      \bfuMh -\divo \Pidiv \bu(t_m)}\,dt
    \\
    &=\Int (\ff (t_m)-\bff(t),\bfuMh -\Pidiv\bu(t_m)  )\, dt\,.
  \end{aligned}
\end{align}
We now discuss and estimate the terms in \eqref{eq:err-fin} separately, to arrive finally
to the estimate~\eqref{eq:error-estimate}.  

First, note that the projection operator $\Pidiv $ has the same properties as the operator
$P_h$ considered in \cite{br-parabolic} and that the solution $\bu$ of \eqref{eq:cont-var}
and the solution treated in \cite{br-parabolic} possess exactly the same regularity. Thus,
the first two terms on the left-hand side can be treated exactly as in
\cite{br-parabolic}. Consequently, \cite[Lemmas 3.7, 3.9]{br-parabolic} yield the following estimates:
\begin{equation}
\begin{aligned}\label{eq:time}
    &\big ({d_t}(\bfuMh-\bfu(t_{m})), \bfuMh-\Pidiv \bfu(t_{m})\big)
    \\
    &\quad\geq\frac{1}{2}d_t\|\bfuMh-\bfu(t_{m})\|^2_2 -c\,
    \frac{h^{2+4/p'}}{\kappa}\Int \|\nabla
    ^2\bu(t)\|_{\frac{6p}{4+p}}^2\,d t
    \\
    &\qquad -c\,\frac{h^{4/p'}}{\kappa}\Big\|\nabla\bfu(t_{m})-\Int
    \nabla \bfu(t)\,dt\Big\|_{\frac{6p}{4+p}}^2\,,
  \end{aligned}
\end{equation}
and
\begin{equation}
\begin{aligned}\label{eq:S}
    &\Int\bighskp{ \bS(\bD\bfuMh)-\bS(\bD\bfu(t))}{\bD\bfuMh -\bD
      \Pidiv \bfu(t_{m})} \,dt
    \\
    &\quad\geq\|\bF(\bD\bfuMh)-\bF(\bD\bu(t_m))\|^{2}_2
    -c\, h^{2}\Int\|\nabla\bF(\bD\bu(t))\|^{2}_2\,dt \\
    & \qquad - c\Int\|\bF(\bD\bu(t))-\bF(\bD\bu(t_{m}))\|^2_2\,dt\,.
  \end{aligned}
\end{equation}

The term with the external force is treated slightly differently
compared to \cite[Lemma 3.10]{br-parabolic}. This is due to the fact
that to apply Gronwall's lemma~\ref{lem:discrete_Gronwall_lemma} we
need an estimate involving the $L^{p}$-norm of the gradient of the
error. To shorten the notation in the following computations we denote
the error for $m=1,\ldots, M$ by
\begin{align*}
  \be_h^m:=\bfuMh-\bu(t_m)\,.
\end{align*}

\begin{lemma}
  \label{lem:f}
  Under the hypotheses of Theorem \ref{thm:theorem-parabolic} we have
\begin{align*}
  &\Bigabs{ \;\Int \big(\ff(t_m) - {\bff}(t),\bfuMh-\Pidiv \bu(t_m)\big )\,dt\,}
  \\
  &\le c\, \Int\| \ff(t_m) - {\bff}(t) \|_2^2\, dt  + c\, \|\bD\bfeMh\|_p\, \Int\|
    \ff(t_m) - {\bff}(t) \|_2\, dt 
    \\
    &\quad +c\,h^{2+4/p'} \Int \|\nabla ^2 \bu(t)\|_{\frac{6p}{p+4}}
    ^2\,dt +c\,h^{4/p'}\Int \|\nabla \bfu(t_{m})-
    \nabla \bfu(t)\|_{\frac{6p}{4+p}}^2\, dt.
\end{align*}
\end{lemma}
\begin{proof}
  Using that 
  \begin{equation*}
    \bfuMh -\Pidiv \bu(t_m) =\bfuMh-\bfu(t_{m})+\bfu(t_{m})-\Pidiv \bu(t_m) =\bfeMh+\bfu(t_{m})-\Pidiv \bu(t_m)\,,
  \end{equation*}
  together with H\"older's inequality,  Young's inequality and the
  embedding $W^{1,p}(\Omega)\vnor L^2(\Omega)$, valid for $p \ge \frac 65$, we get
\begin{align*}
  &\Bigabs{ \;\Int \big(\ff(t_m) - {\bff}(t),\bfuMh-\Pidiv\bu(t_m)\big )\,dt}
  \\
  &\le c\,     \|\bD\bfeMh\|_p\,\Int\|    \ff(t_m) - {\bff}(t) \|_2\, dt  
  \\
  &\quad + c\, \Int\| \ff(t_m) - {\bff}(t) \|_2^2\, dt  + c\,
    \|\bu(t_m)-\Pidiv \bu(t_m)\|_2^2\,.
\end{align*}
The last term was already treated in the proof of
\cite[Lemma~3.7]{br-parabolic}, where it is proved that  
\begin{align*}
  &\|\bu(t_m)-\Pidiv \bu(t_m)\|_2^2
  \\
  &\le c\,h^{2+4/p'} \Int \|\nabla ^2 \bu(t)\|_{\frac{6p}{p+4}}
    ^2\,dt +c\,h^{4/p'}\Int \|\nabla \bfu(t_{m})-
    \nabla \bfu(t)\|_{\frac{6p}{4+p}}^2\, dt,
\end{align*}
which yields the assertion.
\end{proof}
It remains to treat the convective term and the pressure term, which were not present
in~\cite{br-parabolic} and which require a precise estimation. Let us start with the
former one.
\begin{lemma}\label{lem:conv}
  Under the hypotheses of Theorem \ref{thm:theorem-parabolic} we have
  \begin{align}
    \label{eq:conv}
      &\Bigabs{\;\Int b(\bfuMsh,\bfuMh, \bfuMh -\Pidiv \bu(t_m)) -
      b(\bu(t),\bu(t), \bfuMh -\Pidiv \bu(t_m) )\, dt}
      \\
      &\le c\,  \norm{\bD \bfeMh}_{p}\norm{\bD\bfeMsh}_p^{1-\theta}\norm{\bfeMsh}_2^{\theta}
     +c\, 
     \norm{\bD  \bfeMh}_{p}\,\Int \|\nabla \bu(t_m)-\nabla \bu(t) \|_{\frac
       {6p}{4+p}} \, dt \notag
     \\
     &\quad 
   +c\, \norm{\bD  \bfeMh}_{p}\,\Int 
    \norm{\nabla  \bu(t)-\nabla \bu(t_{m-1})}_{\frac {6p}{4+p}}\, dt
    +c\, \norm{\bD  \bfeMh}_{p} \,h\,\Int \big \|\nabla ^2
   \bu(t)\big \|_{\frac {6p}{4+p}}\, dt 
   \,.\notag
  \end{align}
  Moreover, the estimate is also correct if the first term on the
  right-hand side is replaced by
  $c\,\norm{\bD\bfeMh}_p \norm{\bD \bfeMsh}_{p} $.
\end{lemma}

\begin{proof}
Since 
${\Pidiv \bfuMh=\bfuMh}$, we get $\bfuMh-\Pidiv \bfu(t_{m})=\Pidiv\be_{h}^{m}$.
Thus,  we can re-write the integrand in the term to be estimated in \eqref{eq:conv} as follows
 \begin{equation*}
   b(\bfuMsh,\bfuMh,\Pidiv\be_{h}^{m}) -
      b(\bu(t),\bu(t), \Pidiv\be_{h}^{m} )\,.
 \end{equation*}
To the latter we add and subtract, in the order, the terms 
$b( \bfuMsh,\Pidiv \bu(t_m),\Pidiv\be_{h}^{m})$,
$b( \bu(t_{m-1}), \Pidiv \bu(t_m), \Pidiv\be_{h}^{m})$,
$b( \bu(t_{m-1}),\bu(t_m), \Pidiv\be_{h}^{m})$ and 
${b(\bu (t), \bu(t_m), \Pidiv\be_{h}^{m})}$, to get for all 
$m=1,\ldots , M$ and a.e.~$t \in I_m$
\begin{align}
  \label{eq:abc}
  \begin{aligned}
    &\hspace{-1cm}b(\bfuMsh,\bfuMh, \Pidiv\be_{h}^{m}) - b(\bu(t),\bu(t), \Pidiv\be_{h}^{m} )
      \\
      &=b(\bfuMsh,\bfuMh -\Pidiv \bu(t_m) , \Pidiv\be_{h}^{m} )
      \\
      &\quad + b(\bfuMsh-\bu(t_{m-1}),\Pidiv\bu(t_m), \Pidiv\be_{h}^{m}  )
      \\
      &\quad +b( \bu(t_{m-1}), \Pidiv \bu(t_m)-\bu(t_m),\Pidiv\be_{h}^{m} )
      \\
      &\quad + b(\bu(t_{m-1})-\bu(t),\bu(t_m), \Pidiv\be_{h}^{m}  )
      \\
      &\quad + b(\bu (t), \bu(t_m)-\bu(t), \Pidiv\be_{h}^{m})
      \\
      &=: I_1^m(t)+I_2^m(t)+I_3^m(t)+I_4^m(t)+I_5^m(t)\,.
  \end{aligned}
\end{align}
In view of the skew-symmetry of $b(\cdot,\cdot,\cdot)$ we have
$I_1^m(t)=0$, for all $m=1,\ldots , M$ and $t\in I_m$.

Using the definition of $b(\cdot,\cdot,\cdot)$ in \eqref{def:b}, and partial integration
we get for all $m=1,\ldots , M$ and $t\in I_m$, also using H\"older's inequality with
$(\frac{3p}2,\frac{3p}{4p-5}, \frac {3p}{3-p})$ and $(p,\frac{3p}{4p-5}, \frac {3p}{2-p})
$, respectively, the embeddings $W^{1,p}(\Omega) \vnor L^{\frac {3p}{3-p}}(\Omega)$,
$W^{1,\frac {3p}2}(\Omega) \vnor L^{\frac {3p}{2-p}}(\Omega)$, Korn's inequality, the
interpolation of $L^{\frac{3p}{4p-5}}(\Omega)$ between $L^2(\Omega)$ and
$W^{1,p}(\Omega)$, which is possible for $p \in (\frac 85,2]$, the continuity of $\Pidiv $
(cf.~Proposition~\ref{prop:Ph}~(ii)), and
\eqref{eq:I-interpolation-u}, that
\begin{align}
  \label{eq:j2}
     \abs{I_2^m(t)} &\le \frac 12 \bigabs{ \hskp{[\nabla \Pidiv
        \bu(t_m)]\bfeMsh}{\Pidiv \bfeMh }}
    + \frac 12 \bigabs{\hskp{[\nabla \Pidiv \bfeMh
        ]\bfeMsh}{\Pidiv \bu(t_m)}} \notag
  \\
      &\le c\, \norm{\nabla \Pidiv \bu(t_m)}_{\frac {3p}2}
    \norm{\bfeMsh}_{\frac {3p}{4p-5}}
    \norm{\bD \Pidiv \bfeMh}_{p}
    \\
    &\le c\, 
    \norm{\bfeMsh}_2^{\theta}  \norm{\bD\bfeMsh}_p^{1-\theta}
    \norm{\bD  \bfeMh}_{p}\,, \notag 
\end{align}
with $\theta :=\frac{10p-16}{5p-6} \in (0,1] $ for $p \in (\frac 85
,2]$. Using the embedding $W^{1,p}(\Omega) \vnor L^2(\Omega)$ in the
last line we also obtain 
\begin{align}
  \label{eq:j2a}
  \begin{aligned}
    \abs{I_2^m(t)} &\le c\, 
    \norm{\bD\bfeMsh}_p \norm{\bD  \bfeMh}_{p}\,.
  \end{aligned}
\end{align}
Since $\bu(t_{m-1})$ is solenoidal we get, also
using H\"older's inequality with
$(\frac {3p}{5(p-1)},\frac{3p}{2-p}, \frac {3p}{3-p})$, the embeddings
$W^{1,p}(\Omega) \vnor L^{\frac {3p}{3-p}}(\Omega)$,
$W^{1,\frac {3p}2}(\Omega) \vnor L^{\frac
  {3p}{2-p}}(\Omega)$, 
$L^{\frac {6p}{4+p}}(\Omega) \vnor L^{\frac {3p}{5(p-1)}}(\Omega)$, 
the continuity of $\Pidiv $
(cf.~Proposition \ref{prop:Ph} (ii)), and
\eqref{eq:I-interpolation-u}, that
\begin{align}
  \label{eq:j3}
  \begin{aligned}
    \abs{I_3^m(t)} &= \bigabs{ \bighskp{[\nabla \Pidiv \bu(t_m)-\nabla\bu(t_m)]
        \bfu(t_{m-1})}{\Pidiv \bfeMh }} 
    \\
    &\le c\, \norm{\nabla \Pidiv \bu(t_m)-\nabla \bu(t_m)}_{\frac {3p}{5(p-1)}}
    \norm{\nabla \bu(t_m)}_{\frac {3p}2}    \norm{\bD \Pidiv \bfeMh}_{p}
    \\
    &\le c\, 
    \norm{\nabla \Pidiv \bu(t_m)-\nabla \bu(t_m)}_{\frac {6p}{4+p}}
    \norm{\bD  \bfeMh}_{p}\,.
  \end{aligned}
\end{align}
To treat
$\norm{\nabla \Pidiv \bu(t_m)-\nabla\bu(t_m)}_{\frac {6p}{4+p}}$ we
add and subtract $\nabla \Pidiv \big(\Int \bu(\sigma )\, d\sigma \big )$, use
$\Pidiv \big(\Int \bu(\sigma )\, d\sigma \big )= \Int \Pidiv \bu (\sigma )\, d\sigma$, add and subtract
$\nabla \Int\bu(\sigma )\, d\sigma$, use the continuity of $\Pidiv $
(cf.~Proposition \ref{prop:Ph}~(ii)),
$\nabla \Int \bv(\sigma )\, d\sigma = \Int\nabla \bv(\sigma )\, d\sigma$, Fubini's theorem,
the properties of the Bochner integral, Proposition \ref{prop:Ph}~(ii)
and \eqref{eq:regularity} to arrive at
\begin{align}
  \label{eq:j3a}
  \begin{aligned}
    &\norm{\nabla \Pidiv \bu(t_m)-\nabla \bu(t_m))}_{\frac {6p}{4+p}}
    \\
    &\le c\, \big \|\nabla \Pidiv \Int \bu(t_m)-\bu(\sigma )\, d\sigma \big
    \|_{\frac {6p}{4+p}} +c\, \big \|\nabla \Int \Pidiv
    \bu(\sigma )-\bu(\sigma ) \, d\sigma \big \|_{\frac {6p}{4+p}}
    \\
    &\quad +c\, \big \|\nabla \Int \bu(\sigma ) -\bu(t_m) \, d\sigma\big
    \|_{\frac {6p}{4+p}}
    \\
    &\le c\, \big \|\Int \nabla \bu(t_m)-\nabla \bu(\sigma )\, d\sigma \big
    \|_{\frac {6p}{4+p}} +c\, \big \| \Int \nabla \Pidiv
    \bu(\sigma )-\nabla \bu(\sigma ) \, d\sigma \big \|_{\frac {6p}{4+p}}
    \\
    &\le c\, \Int \big \|\nabla \bu(t_m)-\nabla \bu(\sigma )\big
    \|_{\frac {6p}{4+p}}\, d\sigma  +c\, \Int \big \|\nabla \Pidiv
    \bu(\sigma )-\nabla \bu(\sigma ) \big \|_{\frac {6p}{4+p}}\, d\sigma 
    \\
    &\le c\, \Int \big \|\nabla \bu(t_m)-\nabla \bu(\sigma )\big
    \|_{\frac {6p}{4+p}}\, d\sigma  +c\, h\,\Int \big \|\nabla ^2
    \bu(\sigma )\big \|_{\frac {6p}{4+p}}\, d\sigma \,.
  \end{aligned}
\end{align}
Combining \eqref{eq:j3} and \eqref{eq:j3a} we showed that for all
$m=1,\ldots , M$ and $t\in I_m$ there holds 
\begin{align}
  \label{eq:j3b}
  \begin{aligned}
    \abs{I_3^m(t)}&\le c\, 
     \norm{\bD  \bfeMh}_{p}\,\Int \|\nabla \bu(t_m)-\nabla \bu(\sigma ) \|_{\frac
       {6p}{4+p}} \, d\sigma 
     \\
     &\quad  +c\, \norm{\bD  \bfeMh}_{p} \, h\,\Int \big \|\nabla ^2
    \bu(\sigma )\big \|_{\frac {6p}{4+p}}\, d\sigma \,.
  \end{aligned}
\end{align}
Since $\bu(t_{m-1})$ and $\bu(t)$ are solenoidal we get for every  $m=1,\ldots , M$ and
$t \in I_m$, also using H\"older's inequality with
$(\frac{3p}2,\frac{3p}{4p-5}, \frac {3p}{3-p})$, the embedding
$W^{1,p}(\Omega) \vnor L^{\frac {3p}{3-p}}(\Omega)$, Korn's
inequality, the embedding
$W^{1,\frac {6p}{4+p}}(\Omega) \vnor L^{\frac {6p}{4-p}}(\Omega) \vnor L^{\frac {3p}{4p-5}}(\Omega)$,
the continuity of $\Pidiv $ (cf.~Proposition \ref{prop:Ph} (ii)), and
\eqref{eq:I-interpolation-u}, that
\begin{align}
  \label{eq:j4}
  \begin{aligned}
    \abs{I_4^m(t)} &= \bigabs{ \bighskp{[\nabla \bu(t_m)]\big (\bu(t)-
        \bfu(t_{m-1})\big)}{\Pidiv \bfeMh }} 
    \\
    &\le c\, \norm{\nabla \bu(t_m)}_{\frac {3p}{2}}
    \norm{\bu(t)-\bu(t_m)}_{\frac {3p}{4p-5}}    \norm{\bD \Pidiv \bfeMh}_{p}
    \\
    &\le c\, 
    \norm{\nabla  \bu(t)-\nabla \bu(t_{m-1})}_{\frac {6p}{4+p}}
    \norm{\bD  \bfeMh}_{p}\,.
  \end{aligned}
\end{align}
Since $\bu(t)$ is solenoidal we get for every  $m=1,\ldots , M$ and $t \in I_m$, also using
H\"older's inequality with
$(\frac{6p}{4+p},\frac{6p}{7p-10}, \frac {3p}{3-p})$, the embedding
$W^{1,p}(\Omega) \vnor L^{\frac {3p}{3-p}}(\Omega)$, Korn's
inequality, the embedding
$W^{1,\frac {3p}2}(\Omega) \vnor L^{\frac {3p}{2-p}}(\Omega) \vnor L^{\frac {6p}{7p-10}}(\Omega)$,
\eqref{eq:I-interpolation-u}, and the
continuity of $\Pidiv $ (cf.~Proposition \ref{prop:Ph} (ii)), that
\begin{align}
  \label{eq:j5}
  \begin{aligned}
    \abs{I_5^m(t)} &= \bigabs{ \bighskp{[\nabla \bu(t_m)- \nabla \bu(t)]\bu(t)}{\Pidiv \bfeMh }} 
    \\
    &\le c\, \norm{\nabla \bu(t_m)-\nabla \bu(t)}_{\frac {6p}{4+p}}
    \norm{\bu(t)}_{\frac {6p}{7p-10}}    \norm{\bD \Pidiv \bfeMh}_{p}
    \\
    &\le c\, 
    \norm{\nabla  \bu(t_m)-\nabla \bu(t)}_{\frac {6p}{4+p}}
    \norm{\bD  \bfeMh}_{p}\,.
  \end{aligned}
\end{align}
The assertions follow from \eqref{eq:j2} and \eqref{eq:j2a}, resp., \eqref{eq:j3b},
\eqref{eq:j4} and \eqref{eq:j5}. 
\end{proof}

Let us now discuss the last term, namely the one including the
pressure.
\begin{lemma}\label{lem:q}
Under the hypotheses of Theorem \ref{thm:theorem-parabolic} for
  every $\vep >0$ there exists $c_\vep>0$ such that 
  \begin{align}
    \label{eq:q}
    \begin{aligned}
     &\Bigabs{ \,\Int \bighskp{\PiY q (t) - q(t)}{\divo
         \bfuMh -\divo\Pidiv \bu(t_m)}\,dt}
     \\
     &\le \vep \, \norm{\bF (\bD \bfuMh) -\bF(\bD\bu(t_m))}_2^2+ c_\vep\,
     \Int \norm{\bF (\bD \bfu(t_m))
       -\bF(\bD\bu(t))}_2^2 \, dt
     \\
     &\quad + c_\vep\, \Int \norm{\bF (\bD \bfu(t))
        -\bF(\bD\Pidiv \bu(t))}_2^2 \, dt +c_\vep\, h^2\Int \norm{\nabla
        \bF(\bD \bu(t))}_2^2 \, dt
      \\
      &\quad +c_\vep\, h^2 \Int \|\nabla q(t)\|_{p'}^{p'} +
       \norm{\bF(\bD \bu(t))}_2^2 \, dt  + c_\vep\, h^2\,\abs{\Omega}\,\delta^{p'}\,.
    \end{aligned}
  \end{align}
\end{lemma}
\begin{proof}
  The  point-wise application of Young's inequality \eqref{ineq:young}
  with $a=\abs{\bD\bu(t_m,\cdot)}$ to the integrand of the term on the
  left-hand side of \eqref{eq:q} yields for all $m=1,\ldots, M$, 
  a.e.~$t \in I_m$ and every $\vep >0$
  \begin{align}
    \label{eq:q1}
    &\bigabs{\bighskp{\PiY q (t) - q(t)}{\divo
      \bfuMh -\Pidiv \bu(t_m)}}
    \\
    &\le  \frac{\vep}{c_{0}} \int\limits_ \Omega\!
      \phi_{\abs{\bD\bu(t_m)}}(\abs{\bD\bfuMh \!-\!\bD \Pidiv \bu(t_m)})\, dx
      +c_\vep \int\limits_ \Omega \! \big(\phi_{\abs{\bD\bu(t_m)}}\big
      )^*(\abs{\PiY q(t) \!-\! q(t)})\, dx 
      \notag
    \\
    &\le  \vep \, \norm{\bF (\bD \bfuMh)
      -\bF(\bD\bu(t_m))}_2^2+\vep\,  \norm{\bF (\bD \bfu(t_m))
      -\bF(\bD\Pidiv \bu(t_m))}_2^2\notag 
    \\
    &\quad +
      c_\vep \int\limits_ \Omega  \big(\phi_{\abs{\bD\bu(t_m)}}\big
      )^*(\abs{\PiY q(t) - q(t)})\, dx \notag
    \\
    &=: \vep\, J^m_1(t) + \vep\, J^m_2(t) + c_\vep J_3^m(t)\,,\notag 
  \end{align}
  where in the last estimate we added and subtracted $\bD\bu(t_m)$ and
  used that for all $\bP, \bQ \in \setR^{3\times 3}_{\sym}$ there
  holds $ \phi_{\abs{\bfP}}(\abs{\bfP - \bfQ}) \le  c_{0}
  \bigabs{ \bfF(\bfP) - \bfF(\bfQ)}^2$
  (cf.~Proposition~\ref{lem:hammer}~(i)). To treat $J_2^m(t)$ we add and
  subtract $\bF(\bD\bu(t))$ and $\bF(\bD\Pidiv \bu(t))$
  \begin{align}
    \label{eq:q2}
    \begin{aligned}
      \abs{J_2^m(t)}&\le c\, \norm{\bF (\bD \bfu(t_m))
        -\bF(\bD\bu(t))}_2^2 + \norm{\bF (\bD \bfu(t)) -\bF(\bD\Pidiv
        \bu(t))}_2^2
      \\
      &\quad + \norm{\bF (\bD \Pidiv \bfu(t)) -\bF(\bD\Pidiv
        \bu(t_m))}_2^2\,.
    \end{aligned}
  \end{align}
  To treat the last term on the right-hand side of \eqref{eq:q2} we use Proposition
  \ref{lem:hammer} and Proposition~\ref{prop:Ph}~(iv) to get
  \begin{align}
    \label{eq:q3}
      &\norm{\bF (\bD \Pidiv \bfu(t)) -\bF(\bD\Pidiv \bu(t_m))}_2^2 
      \\
      & \le c\,\int\limits _\Omega \phi_{\abs{\bD\Pidiv
          \bu(t)}}(\abs{\bD\Pidiv \bfu(t)-\bD \Pidiv \bu(t_m)})\, dx \notag 
      \\
      & \le c\,\int\limits _\Omega \phi_{\abs{\bD
          \bu(t)}}(\abs{\bD\Pidiv \bfu(t)-\bD \Pidiv \bu(t_m)})\, dx
       + c\, \norm{\bF (\bD \Pidiv \bfu(t)) -\bF(\bD \bu(t))}_2^2 \notag 
      \\
      &\le c\, h^2\, \norm{\nabla \bF(\bD \bu(t))}_2^2 + c\, \norm{\bF
        (\bD \bfu(t)) -\bF(\bD \bu(t_m))}_2^2 \notag 
      \\
      &\quad   + c\, \norm{\bF (\bD \Pidiv \bfu(t)) -\bF(\bD
        \bu(t))}_2^2\,. \notag 
  \end{align}
  To treat the term $J_3^m(t)$ we note that $\Omega = \bigcup _{K \in
    \mathcal T_h} K$,  use Proposition~\ref{lem:hammer} (ii), (iii)
  and $K \subseteq S_K$ to arrive at 
  \begin{align}
    \label{eq:q-err}
    \begin{aligned}
      \abs{J_3^m(t)}&\le c \sum_{K \in \mathcal T_h} \int\limits_K \big(\phi_{\abs{\mean{\bD\bu(t_m)}_{S_K}}}\big
      )^*(\abs{\PiY q(t) - q(t)})\, dx
      \\
      &\quad + c \sum_{K \in \mathcal T_h}\;
      \int\limits_{S_K} \abs{\bF(\bD \bu(t_m))-\mean{\bF(\bD
          \bu(t_m))}_{S_K}}^2\, dx
      \\
      &=: c \sum_{K \in \mathcal T_h} A_K^m(t) + c \sum_{K \in \mathcal T_h}
      B_K^m(t)\,.
    \end{aligned}
  \end{align}
  Using Proposition \ref{pro:PY}, again Proposition \ref{lem:hammer}
  (ii), (iii) and $\big(\phi_{\abs{\bD\bu(t) }}\big)^*(h \abs{\sigma} ) \le c\,
  h^2 \big( \abs{\sigma}^{p'} +\delta^{p'}+ \abs{\bF(\bD\bu (t))}^2 \big ) $,
  valid for $p \le 2$,  and $h_{K}\leq h$ yields
  \begin{align}
    \label{eq:q-err1}
    \begin{aligned}
      \abs{A_K^m(t)}&\le \int\limits_{S_K}
      \big(\phi_{\abs{\mean{\bD\bu(t_m)}_{S_K}}}\big
      )^*(h_K\abs{\nabla q(t)})\, dx
      \\
      &\le c \int\limits_{S_K} \big(\phi_{\abs{\bD\bu(t_m)}}\big
      )^*(h_K\abs{\nabla q(t)})\, dx + c\,B_K^m(t)
      \\
      &\le c \int\limits_{S_K} \big(\phi_{\abs{\bD\bu(t)}}\big
      )^*(h_K\abs{\nabla q(t)})\, dx
      \\
      &\quad + \int\limits_{S_K} \abs{\bF(\bD \bu(t_m))-\bF(\bD
        \bu(t))}^2\, dx + c\,B_K^m(t)\,.
      \\
      &\le c \, h^2 \int\limits_{S_K} \abs{\nabla q(t)}^{p'} +
      \delta^{p'} +c\,   \abs{\bF(\bD \bu(t))}^2\, dx
      \\
      &\quad + \int\limits_{S_K} \abs{\bF(\bD \bu(t_m))-\bF(\bD
        \bu(t))}^2\, dx + c\,B_K^m(t)\,.
    \end{aligned}
  \end{align}
  Adding and subtracting appropriate terms, using Proposition
  \ref{lem:hammer}~(iii), the properties of
  the mean value, $\abs{S_K}\sim \abs{K}$ and Poincar\'e's inequality we get
  \begin{align}
    \label{eq:q-err2}
      \abs{B_K^m(t)}&\le c \int\limits_{S_K} \!\abs{\bF(\bD \bu(t_m))-\bF(\bD
        \bu(t))}^2\, dx
    + c \int\limits_{S_K} \!\abs{\bF(\bD \bu(t))-\mean{\bF(\bD
          \bu(t))}_{S_K}}^2\, dx\notag 
      \\
      &\quad + c \int\limits_{S_K} \abs{\mean{\bF(\bD
          \bu(t))}_{S_K} -\mean{\bF(\bD\bu(t_m))}_{S_K}}^2\, dx
      \\
      &\le c \int\limits_{S_K} \abs{\bF(\bD \bu(t_m))-\bF(\bD
        \bu(t))}^2\, dx  + c\, h^2\int\limits_{S_K} \abs{\nabla
        \bF(\bD \bu(t))}^2\, dx \,.\notag 
  \end{align}
  The assertion follows from \eqref{eq:q1}--\eqref{eq:q-err2} and the
  properties of the triangulation.
\end{proof}

Collecting all estimates,  we are ready to prove the main result of this paper.
\begin{proof}[Proof of Theorem~\ref{thm:theorem-parabolic}]
  From estimates~\eqref{eq:time} and \eqref{eq:S}, Lemma~\ref{lem:f},
  Lemma~\ref{lem:conv}, and Lemma~\ref{lem:q} we get, choosing $\vep >0$ sufficiently
  small to absorb the term $\vep \, \norm{\bF (\bD \bfuMh) -\bF(\bD\bu(t_m))}_2^2$, and
  using the properties of the retarded time averages, that for all $m=1,\dots,M$ and
  $0<\kappa\le 1 $
\begin{align}
  \label{eq:error-estimate}
    &{d_t}\|\bfuMh-\bfu(t_m)\|_2^2+c\,\|\bF(\bD\bfuMh)-\bF(\bD\bfu(t_{m}))\|^2_2
    \\
    &\leq c\, h^{2}\Int\|\nabla\bF(\bD\bu(t) )\|_2^{2}\,dt + c\, \frac
    {h^{2+4/p'}}\kappa \Int \| \nabla ^2\bu(t)\|_{\frac{6p}{4+p}}
    ^2\,dt +c\, h^2 \Int \|\nabla q(t)\|_{p'}^{p'} \, dt \notag
    \\
    &\quad+ c\Int\|\bF(\bD\bu(t))-\bF(\bD\bu(t_{m}))\|^2_2\,dt +
    c\, \frac{h^{4/p'}}{\kappa}\Int \|\nabla \bfu(t_{m})-
    \nabla\bfu(t)\|_{\frac{6p}{4+p}}^2 \,dt \notag
    \\
    &\quad + c\, \Int \norm{\bF (\bD \bfu(t)) -\bF(\bD\Pidiv
      \bu(t))}_2^2 \, dt +c\,  h^2 \Int \norm{\bF(\bD \bu(t))}_2^2 \, dt\notag
    \\
    &\quad + c\,
    h^2\,\abs{\Omega}\,\delta^{p'}+ c\, \Int
    \|\ff(t_{m})-\ff(t)\|_{2}^2 \,dt + c\, \|\bD\bfeMh\|_p \,\Int\|    \ff(t_m) - {\bff}(t) \|_2\, dt\,
    \notag
    \\
    &\quad +c\, \|\bD\bfeMh\|_p   \norm{\bD\bfeMsh}_p^{1-\theta}\norm{\bfeMsh}_2^{\theta}
     +c\, 
      \norm{\bD \bfeMh}_{p}\,\Int \|\nabla \bu(t_m)-\nabla \bu(t) \|_{\frac
       {6p}{4+p}} \, dt \notag
     \\
     &\quad 
   +c\,  \norm{\bD \bfeMh}_{p} \,h\,\Int \big \|\nabla ^2
   \bu(t)\big \|_{\frac {6p}{4+p}}\, dt 
   +c\,  \norm{\bD \bfeMh}_{p}\,\Int 
    \norm{\nabla  \bu(t)-\nabla \bu(t_{m-1})}_{\frac {6p}{4+p}}\, dt\,. \notag
\end{align}
Moreover, the estimate \eqref{eq:error-estimate} is also correct if
$c\, \norm{\bD \bfeMh}_{p} \norm{\bD\bfeMsh}_p^{1-\theta}
\norm{\bfeMsh}_2^{\theta}$ is replaced by
$c\,\norm{\bD\bfeMh}_p \norm{\bD \bfeMsh}_{p} $.  To use Lemma
\ref{lem:discrete_Gronwall_lemma} we observe that by
Lemma~\ref{july_1} and \eqref{eq:I-interpolation-u} we have with
${\lambda:=\delta+\norm{\bfD \bfu}_{C(\overline {I}; L^p(\Omega))}}$
  \begin{align*}
    \norm{\bfF(\bfD\bfuMh)-\bfF(\bfD\bfu(t_m))}^2_2 &\geq c
     \,\big (\lambda+    \norm{\bfD \bfuMh - \bfD
    \bfu(t_m)}_p\big )^{p -2} \, \norm{ \bfD \bfuMh - \bfD\bfu(t_m)}^2_p
    \\
    &=c \,\big (\lambda+    \norm{\bfD \bfeMh }_p\big )^{p -2} \, \norm{ \bfD \bfeMh }^2_p\,.
  \end{align*}
Thus, the left-hand side of \eqref{eq:error-estimate}  is larger or equal than
  \begin{align*}
    &{d_t}\|\bfeMh\|_2^2 +c \,\big (\lambda+    \norm{\bfD \bfeMh }_p\big )^{p -2} \, \norm{ \bfD \bfeMh }^2_p\,,
  \end{align*}
 and \eqref{eq:error-estimate} is now 
  written in the form needed for the application of Lemma
  \ref{lem:discrete_Gronwall_lemma}.
To do so, we set
\begin{align*}
  a_m(h)&:= \|\bfeMh\|_2\,, \quad b_m(h):=\norm{\bfD \bfeMh }_p\,, \quad
  r_m(h,\kappa):=h\,\Int \big \|\nabla ^2 \bu(t)\big \|_{\frac {6p}{4+p}}\,  dt\,,
\\
  \rho_m(h,\kappa)=\rho_m&:= \Int \|\nabla \bu(t_m)-\nabla \bu(t) \|_{\frac {6p}{4+p}}
  \, dt + \Int \norm{\nabla \bu(t)-\nabla \bu(t_{m-1})}_{\frac
          {6p}{4+p}}\, dt\,
 \\
  &\quad +\Int\|    \ff(t_m) - {\bff}(t) \|_2\, dt\,,
\\
  \big(s_m(h,\kappa)\big)^2&:=h^{2}\Int\|\nabla\bF(\bD\bu(t) )\|_2^{2}\,dt + \frac
    {h^{2+4/p'}}\kappa \Int \| \nabla ^2\bu(t)\|_{\frac{6p}{4+p}}
    ^2\,dt
    \\
    &\quad +\Int \norm{\bF (\bD \bfu(t)) -\bF(\bD\Pidiv
      \bu(t))}_2^2 \, dt +h^2 \Int \|\nabla q(t)\|_{p'}^{p'} \, dt
  \\
  &\quad + h^2 \Int \norm{\bF(\bD \bu(t))}_2^2 \, dt+ 
    h^2\,\abs{\Omega}\,\delta^{p'}\,,
  \\
    \big(\sigma_m(h,\kappa)\big)^2&:= \frac{h^{4/p'}}{\kappa}\Int \|\nabla \bfu(t_{m})-
    \nabla\bfu(t)\|_{\frac{6p}{4+p}}^2 \,dt  +\Int
            \|\ff(t_{m})-\ff(t)\|_{2}^2 \,dt
  \\
  &\quad +\Int\|\bF(\bD\bu(t))-\bF(\bD\bu(t_{m}))\|^2_2\,dt\,.
\end{align*}
Let us verify that these quantities fulfil assumption \eqref{eq:r_m}.  First, we observe
that, by using the regularity \eqref{eq:regularity} on the velocity $\bu$, it holds
  \begin{equation}\label{est:r}
    \ksumo r_m^2(h,\kappa) \le c\, h^{2}\int\limits_{0}^{T} 
    \|\nabla ^2 \bu(t)\|_{\frac{6p}{p+4}}^2\,dt \le c\, h^2\,. 
  \end{equation}
  Second, by using the regularity \eqref{eq:reg-ass}, the
  condition~\eqref{eq:compatbility} and the regularity
  \eqref{eq:regularity}, Proposition \ref{lem:hammer} (i) and the
  regularity \eqref{eq:reg-ass}, and again the
 regularity \eqref{eq:reg-ass}, we obtain
  \begin{align}\label{est:s}
    \begin{aligned}
      \ksumo s_m^2(h,\kappa) &\le c\, h^{2}\int\limits_{0}^{T} \|\nabla \bF(\bD
      \bu(t))\|_2^2\,dt + h^2\frac {h^{4/p'}}\kappa \int\limits_0^T\|
      \nabla ^2\bu(t)\|_{\frac{6p}{4+p}} ^2\,dt
      \\
      &\quad +h^2 \int\limits_0^T\|\nabla q(t)\|_{p'}^{p'} \, dt + h^2
      \int\limits_0^T \norm{\bF(\bD \bu(t))}_2^2 \, dt+
      h^2\,\abs{\Omega}\,\delta^{p'}
      \\
      & \le c\, h^2\,.
    \end{aligned}
  \end{align}
  Third, using H\"older's inequality, several times Lemma~\ref{lem:Bochner-lemma}, the
  regularity \eqref{eq:regularity}, and the assumption on the regularity of $\ff$ we get
  (since $p\leq2)$
    \begin{align}\label{est:rho}
      \begin{aligned}
        \ksumo \rho_m^2 &\le c\, \kappa^{2}\int\limits_{0}^{T} \|
        \partial _t\nabla \bu(t)\|_{\frac{6p}{p+4}}^2\,dt + c\,
        \kappa^{2}\int\limits_{0}^{T} \| \partial
        _t\ff(t)\|_{2}^2\,dt
        \\
        & \le c\, \kappa^2\,.
      \end{aligned}
  \end{align}
  Next, by using Lemma~\ref{lem:Bochner-lemma}, the
  condition~\eqref{eq:compatbility}, the regularity~\eqref{eq:regularity}, as well as the
  regularity~\eqref{eq:regularity}, 
  and the assumption on the regularity of $\ff$  we have
 \begin{align}\label{est:sig}
   \begin{aligned}
     \ksumo \sigma_m^2(h,\kappa) &\le c\, \kappa^2\frac {h^{4/p'}}\kappa
     \int\limits_0^T\|\partial _t \nabla \bu(t)\|_{\frac{6p}{4+p}}
     ^2\,dt + c\, \kappa^{2}\int\limits_{0}^{T} \| \partial
     _t\ff(t)\|_{2}^2\,dt
     \\
     &\quad + c\, \kappa^{2}\int\limits_{0}^{T} \|\partial_t \bF(\bD
     \bu(t))\|_2^2\,dt
     \\
     & \le c\, \kappa^2\,.
   \end{aligned}
 \end{align}
 Finally, since $\bu_h^0=\Pidiv \bu_0$, the regularity of $\bu_0$ and
 Proposition \ref{prop:Ph} yield
 \begin{align*}
   a_0(h)&=\norm{\bu_0 -\Pidiv \bu _0}_2 \le c\, h\, \norm{\nabla \bu
   _0}_2 \le c\, h\,,
   \\
   b_0(h)&=\norm{\bD \bu_0 -\bD \Pidiv \bu _0}_p \le c\, h\, \norm{\nabla ^2\bu
   _0}_p \le c\, h\,.
 \end{align*}
 Consequently Lemma \ref{lem:discrete_Gronwall_lemma} yields for
 sufficiently small $\kappa$ and $h$ that
   \begin{align*}
    \max_{1\leq m\leq M} \|\bfeMh \|_2^2+{\gamma_1(1+\Lambda)^{p-2}}
    \ksum \|\bD \bfeMh \|_p^2 &\leq c\,\big (h^{2} +\kappa^2\big )\,, 
    \\
    \max_{1\leq m\leq M} \|\bD \bfeMh \|_p&\leq 1\,.
   \end{align*}
   Using this, the estimates \eqref{est:r}--\eqref{est:sig} and
   $\norm{\bfeMsh}_2^{\theta} \norm{\bD\bfeMsh}_p^{1-\theta} \norm{\bD
     \bfeMh}_{p} \le \linebreak c\, \norm{\bD\bfeMsh}_p^2 +c\,\norm{\bD
  \bfeMh}_{p}^2 $, one easily checks that all terms on the right-hand
   side of \eqref{eq:error-estimate} are bounded, after multiplication
   by $\kappa$ and summation over $m=1,\ldots, M$, by
   $c\,(h^{2} +\kappa^2)$. Thus, we proved
   \begin{equation*}
  \max_{m=1,\ldots,M}
  \|\bfuMh-\bfu(t_{m})\|_2^2+\ksumo\|\bF(\bD\bfuMh)-\bF(\bD\bfu(t_{m}))\|_2^2\leq
  c\,(h^2+\kappa^2),
\end{equation*}
which is the assertion of Theorem \ref{thm:theorem-parabolic}. 
\end{proof}

\begin{proof}[Proof of Corollary \ref{cor:main}]
  This corollary is proved in the same way as Theorem
  \ref{thm:theorem-parabolic} with the only difference that in all
  places where \eqref{eq:I-interpolation-u} is used we use
  \eqref{ass:reg2} instead.

  Let us illustrate that on the estimate of $I_2^m(t)$ in
  \eqref{eq:j2}.  Using the definition of $b(\cdot,\cdot,\cdot)$ in
  \eqref{def:b}, and partial integration we get for all
  $m=1,\ldots , M$ and $t\in I_m$, also using H\"older's inequality
  with $(r,\frac{3pr}{3pr-2p-3}, \frac {3p}{3-p})$ and
  $(p,\frac{p}{p-1}, \infty) $ for $r \in (3,6(p-1))$, respectively,
  the embeddings $W^{1,p}(\Omega) \vnor L^{\frac {3p}{3-p}}(\Omega)$,
  $W^{1,r}(\Omega) \vnor L^{\infty}(\Omega)$, Korn's inequality, the
  embedding $L^{\frac {p}{p-1}}(\Omega) \vnor L^{\frac
    {3pr}{3pr-2p-3}}(\Omega)$, the
  interpolation of $L^{\frac{p}{p-1}}(\Omega)$ between
  $L^2(\Omega)$ and $W^{1,p}(\Omega)$, which is possible for
  $p \in (\frac 32,\frac 85]$, the continuity of $\Pidiv $
  (cf.~Proposition~\ref{prop:Ph}~(ii)), and \eqref{ass:reg2}
\begin{align*}
     \abs{I_2^m(t)} &\le \frac 12 \bigabs{ \hskp{[\nabla \Pidiv
        \bu(t_m)]\bfeMsh}{\Pidiv \bfeMh }}
    + \frac 12 \bigabs{\hskp{[\nabla \Pidiv \bfeMh
        ]\bfeMsh}{\Pidiv \bu(t_m)}}\notag
    \\
    &\le c\, \norm{\nabla \Pidiv \bu(t_m)}_{r} \big ( 
    \norm{\bfeMsh}_{\frac {3pr}{3pr-2p-3}}+    \norm{\bfeMsh}_{\frac
      {p}{p-1}}\big )
    \norm{\bD \Pidiv \bfeMh}_{p}
    \\
    &\le c\, 
    \norm{\bfeMsh}_2^{\theta}  \norm{\bD\bfeMsh}_p^{1-\theta}
    \norm{\bD  \bfeMh}_{p}\,, \notag 
\end{align*}
with $\theta :=\frac{8p}{5p-6} \in (0,1] $ for $p \in (\frac
32,\frac 85]$. Using the embedding $W^{1,p}(\Omega) \vnor L^2(\Omega)$ in the
last line we also obtain 
\begin{align*}
    \abs{I_2^m(t)} &\le c\, 
    \norm{\bD\bfeMsh}_p \norm{\bD  \bfeMh}_{p}\,.
\end{align*}
Similar adaptations apply to the treatment of the other terms stemming from the
convective term.  This proves the assertion.
\end{proof}
Due to the presence of the term 
$\bigabs{\hskp{[\nabla \Pidiv \bfeMh ]\bfeMsh}{\Pidiv \bu(t_m)}}$ an
extension of the validity of the error estimate for 
$p\le \frac{3}{2}$ with the present technique is impossible, even with further regularity assumptions on $\bu$.

\section*{Acknowledgments}
The research of Luigi C. Berselli that led to the present paper was partially supported
by a grant of the group GNAMPA of INdAM and by the project of the University of Pisa
within the grant PRA$\_{}2018\_{}52$~UNIPI \textit{Energy and regularity: New techniques
  for classical PDE problems.} 

\def\cprime{$'$} \def\cprime{$'$} \def\cprime{$'$}

\end{document}